\documentclass{amsart}
\pdfoutput=1
\usepackage{amsmath}
\usepackage{amssymb}
\usepackage{amsthm}
\usepackage{bm}
\usepackage{accents}
\usepackage{mathtools}
\usepackage{tikz}
\usetikzlibrary{calc}
\usetikzlibrary{decorations.pathmorphing,shapes}
\usetikzlibrary{automata,positioning}
\usepackage{tikz-cd}
\usepackage{forest}
\usepackage{braket} 
\usepackage{listings}
\usepackage{mdframed}
\usepackage{verbatim}
\usepackage{physics2}
\usephysicsmodule{ab,ab.legacy,diagmat,xmat}
\usepackage{derivative}
\usepackage{fixdif}
\usepackage{stmaryrd}
\usepackage[mathcal]{eucal}
\usepackage{stackengine} 

\usepackage{microtype}

\usepackage{algorithmicx}
\usepackage{algpseudocode}
\usepackage{algorithm}

\usepackage[english]{babel} 
\usepackage[backend=biber,style=alphabetic,maxalphanames=4,maxnames=5,hyperref]{biblatex}
\usepackage{xurl}
\usepackage[bookmarks, colorlinks, breaklinks]{hyperref} 
\hypersetup{linkcolor=blue,citecolor=magenta,filecolor=black,urlcolor=blue}
\usepackage{cleveref}
\usepackage{graphicx}
\graphicspath{{./}}

\usepackage{float}
\usepackage{booktabs}
\usepackage[shortlabels]{enumitem}
\setitemize{noitemsep}
\usepackage{csquotes}
\newlist{mydescription}{description}{1}
\setlist[mydescription]{style=nextline,
                        font=\bfseries,
                        labelindent=1cm, 
                        leftmargin =2cm,
                        rightmargin=1cm,
                        topsep     =1ex
                       }

\usepackage{lipsum}

\DeclarePairedDelimiter{\floor}{\lfloor}{\rfloor}
\DeclarePairedDelimiter{\ceil}{\lceil}{\rceil}

\newtheorem{thm}{Theorem}[section]

\newtheorem{prop}[thm]{Proposition}
\newtheorem{lem}[thm]{Lemma}

\newtheorem{claim}[thm]{Claim}

\theoremstyle{definition}
\newtheorem{defn}[thm]{Definition}

\newtheorem{notn}[thm]{Notation}

\newtheorem{conv}[thm]{Convention}

\theoremstyle{remark}
\newtheorem{rmk}[thm]{Remark}

\newtheorem{warn}[thm]{Warning}

\theoremstyle{plain}
\newtheorem*{thm*}{Theorem}
\newtheorem*{prop*}{Proposition}
\newtheorem*{lem*}{Lemma}
\newtheorem*{cor*}{Corollary}
\newtheorem*{conj*}{Conjecture}

\theoremstyle{definition}
\newtheorem*{defn*}{Definition}
\newtheorem*{exer*}{Exercise}
\newtheorem*{defns*}{Definitions}
\newtheorem*{con*}{Construction}
\newtheorem*{exm*}{Example}
\newtheorem*{exms*}{Examples}
\newtheorem*{notn*}{Notation}
\newtheorem*{notns*}{Notations}
\newtheorem*{addm*}{Addendum}

\theoremstyle{remark}
\newtheorem*{rmk*}{Remark}

\newcommand{\A}{\mathbb{A}}
\newcommand{\G}{\mathbb{G}}

\newcommand{\R}{\mathbb{R}}
\newcommand{\C}{\mathbb{C}}
\newcommand{\E}{\mathbb{E}}
\newcommand{\Z}{\mathbb{Z}}
\newcommand{\Q}{\mathbb{Q}}
\newcommand{\T}{\mathbf{T}}

\renewcommand{\L}{\mathbb{L}}
\renewcommand{\P}{\mathbb{P}}
\newcommand{\ba}{\mathbf{a}}
\newcommand{\bd}{\mathbf{d}}

\newcommand{\ep}{\varepsilon}

\newcommand{\mc}[1]{\mathcal{#1}}
\newcommand{\mf}[1]{\mathfrak{#1}}
\newcommand{\mbf}[1]{\mathbf{#1}}
\newcommand{\mr}[1]{\mathrm{#1}}
\newcommand{\on}[1]{\operatorname{#1}}

\newcommand{\ol}[1]{\overline{#1}}
\newcommand{\ul}[1]{\underline{#1}}
\newcommand{\wt}[1]{\widetilde{#1}}
\newcommand{\wh}[1]{\widehat{#1}}

\newcommand{\bmu}{\bm{\mu}}
\newcommand{\pre}{\mr{pre}}
\newcommand{\vir}{\mr{vir}}
\newcommand{\fl}{\mr{fl}}
\newcommand{\mv}{\mr{mv}}

\DeclareMathOperator{\Hom}{Hom}

\DeclareMathOperator{\Ext}{Ext}

\DeclareMathOperator{\Aut}{Aut}

\DeclareMathOperator{\Spec}{Spec}

\DeclareMathOperator{\Bl}{Bl}

\title{MSP theory for smooth Calabi-Yau threefolds in weighted $\P^4$}
\author{Patrick Lei}
\date{\today}

\begin{document}
    
\begin{abstract}
  We develop the theory of $N$--mixed-spin-$P$ fields for Fermat-type hypersurfaces in $\P(1,1,1,1,2), \P(1,1,1,1,4)$, and $\P(1,1,1,2,5)$, following the theory developed in~\cite{nmsp} for the quintic threefold.
\end{abstract}

\maketitle

\tableofcontents

\section{Introduction}
\label{sec:introduction}

Let $\ba = (1,1,1,1,2)$, $(1,1,1,1,4)$, or $(1,1,1,2,5)$. Let $k = \sum_{i=1}^5 a_i$ and denote the degree $k$ hypersurface in $\P(\ba)$ by $Z$. We will consider in particular the \textit{Fermat-type} hypersurface, which is given by the equation
\[ W \coloneqq \sum_{i=1}^{5} x_i^{b_i} = 0, \qquad \text{where} \qquad b_i \coloneqq \frac{k}{a_i}. \]
By the adjunction formula, we see that $Z$ is a Calabi-Yau threefold. This is the first of a series of papers which studies the all-genus Gromov-Witten theory of $Z$.

The Gromov-Witten theory of $Z$ is an example of a gauged linear sigma model (GLSM), which was introduced by Witten~\cite{mirrorandglsm}. A mathematical construction of GLSMs was completed by Favero-Kim~\cite{faverokim} using the moduli spaces constructed by Fan-Jarvis-Ruan~\cite{glsm}. As a theory based on geometric invariant theory, a GLSM may have several phases corresponding to different GIT stability chambers. For example, the Gromov-Witten theory of $Z$ has another phase, which is the FJRW theory of the quasi-homogeneous polynomial $W$.
Physically, these correspond to different points on the so-called K\"ahler moduli space, so we would also like a mathematical construction of something resembling a global theory of a GLSM.

A geometric construction for varying the K\"ahler parameters and connecting different phases of a GLSM is the theory of Mixed Spin $P$-fields, developed by Chang-Li-Li-Liu~\cite{mspfermat, msp2} for the quintic based on the master space construction of Thaddeus~\cite{gitflips}. Using this method, Chang-Guo-Li-Zhou~\cite{genusonemsp} reproduced earlier results on the Gromov-Witten theory of the quintic threefold and Guo-Ross~\cite{g1mslg,g1lgcy} proved the genus one mirror theorem for the FJRW theory of $\sum_{i=1}^5 x_i^5$ and the LG/CY correspondence for the quintic in genus one.

In order to simplify the computations in large genus, a key insight of Guo was to duplicate the field which corresponds to deformation of the K\"ahler parameter. This idea was realized in~\cite{nmsp} as $N$-Mixed Spin $P$-field theory\footnote{In this paper, we will refer to the theory with an arbitrary $N$ as simply MSP theory and treat it as depending on a parameter $N$.} and was used by Chang-Guo-Li~\cite{nmsp2, nmsp3} to prove the Yamaguchi-Yau conjectures and the BCOV conjecture for the quintic threefold. 

For arbitrary GLSMs without strictly semistable points, a theory of MSP fields was developed by Chang-Guo-Li-Li-Zhou~\cite{generalmsp}.\footnote{In fact, their construction allows them to take GIT quotients of arbitrary affine schemes, but they do not consider a superpotential.} In particular, this theory applies to our examples. The goal of this paper is to explore the property of the MSP moduli space for the Gromov-Witten theory of $Z$. The first property, proved in~\cite{generalmsp}, is that the moduli space $\mc{W}_{g,\gamma,\bd}$ of MSP fields with discrete data $(g,\gamma, \bd)$ is separated and of finite type. In~\Cref{sec:nmsp}, we will prove that our moduli space is an example of the construction of Chang-Guo-Li-Li-Zhou.

We will then equip $\mc{W}_{g,\gamma,\bd}$ with an action of $\T = (\C^{\times})^N$ and give $\mc{W}_{g,\gamma,\bd}$ a virtual cycle using the cosection localization construction of Kiem-Li~\cite{cosection}, where the cosection is constructed using the polynomial $W$.
\begin{thm}
    The stack $\mc{W}_{g,\gamma,\bd}$ has a virtual fundamental class supported on a proper closed substack $\mc{W}_{g,\gamma,\bd}^- \subset \mc{W}_{g,\gamma,\bd}$.
\end{thm}

A standard technique to compute enumerative invariants is virtual localization. In~\Cref{sec:T fixed loci and graphs}, we will give a combinatorial description of the fixed loci $\mc{W}_{g,\gamma,\bd}^{\T}$ in terms of localization graphs. In particular, the vertices of any localization graph $\Theta$ may be classified by discrete data called the \textit{level}, which may be $0$, $1$, or $\infty$.
\begin{itemize}
    \item The level $0$ vertices correspond to the Gromov-Witten theory of $Z$;
    \item The level $1$ vertices correspond to Hodge integrals;
    \item The level $\infty$ vertices correspond to the FJRW theory of $W$.
\end{itemize}
This is a concrete realization of the principle that MSP theory connects different phases of a GLSM.

There are some localization graphs $\Theta$ which are very difficult to deal with, which we will call \textit{irregular}. This includes all $\Theta$ which contain an edge between a vertex at level $0$ and a vertex at level $\infty$. Fortunately, in the case of the quintic, it was discovered by Chang-Li~\cite{msp3} that these graphs do not contribute to MSP invariants. In~\Cref{sec:irrvan}, we prove the analogous result in our setting:
\begin{thm}
    Let $\Theta$ be irregular. Then if $\mc{W}_{(\Theta)} \subset \mc{W}_{g,\gamma,\bd}^{\T}$ is the fixed component corresponding to $\Theta$, we have
    \[ [W_{(\Theta)}]^{\vir} = 0. \]
\end{thm}

Finally, in~\Cref{sec:Virtual localization}, we will state and prove the virtual localization formula for $\mc{W}_{g,\gamma,\bd}$. In~\Cref{sec:properties}, we will provide an alternative proof of the fact that $\mc{W}_{g,\gamma,\bd}$ is separated and that $\mc{W}_{g,\gamma,\bd}^-$ is proper and correct some slight inaccuracies in~\cite{mspfermat}.
In subsequent work~\cite{polynomiality,bcovme}, we will prove the Yamaguchi-Yau polynomial structure conjecture and the BCOV conjecture for the Gromov-Witten theory of $Z$.

\begin{conv}
We define algebraic stacks and Deligne-Mumford stacks following~\cite{champs}. In particular, the diagonal of an algebraic stack is always separated and quasi-compact.
\end{conv}

\subsection*{Acknowledgements}%
\label{sub:Acknowledgements}

The author would like to thank Chiu-Chu Melissa Liu for proposing this project and explaining virtual localization as well as Huai-Liang Chang, Jun Li, and Yang Zhou for helpful discussions about irregular vanishing. The author would also like to thank Jun Li and Yang Zhou again for their hospitality during the author's visit to the Shanghai Center for Mathematical Sciences in May 2024, where some of this work was completed.

\section{MSP fields}
\label{sec:nmsp}

Let $\bm{ \mu }_k \leq \C^{\times}$ be the subgroup of $k$-th roots of unity. We define
\[ \bm{\mu}_k^{\mr{br}} \coloneqq \bm{\mu}_k \cup \bigcup_{i} \bmu_{a_i} \times \ab\{\rho, \varphi\}  \qquad \bm{\mu}_k^{\mr{na}} \coloneqq \bm{\mu}_k^{\mr{br}} \setminus\ab(\bigcup_i \bmu_{a_i}). \]
Our discrete data will consist of a tuple $(g,\gamma, \mbf{d})$, where
\[ g \in \Z_{\geq 0}, \qquad \gamma \in (\bm{\mu}_k^{\mr{br}})^{\ell}, \qquad \mbf{d} = (d_0, d_{\infty}) \in \Q^2. \]
If all $\gamma_i \in \bmu_k^{\mr{na}}$, then we call $(g, \gamma,  \mbf{d})$ \textit{narrow}. Otherwise, we call it \textit{broad}.

\begin{rmk}
    While all results can be written in greater generality, from now on we will only consider $\gamma \in \ab\{ (1,\rho) \} \cup \bmu_k \setminus \ab(\bigcup_i \bmu_{a_i})$, which will correspond to either GW insertions or narrow FJRW insertions. Recall that all broad FJRW invariants vanish~\cite[Theorem 2.1]{wittentopchern}.
\end{rmk}

\begin{defn}\label{defn:nmsp}
  Let $S$ be a scheme. An $S$-family of \textit{weighted MSP fields} of type $(g,\gamma, \mbf{d})$ is a tuple
  \[ \xi = (\mc{C}, \Sigma^{\mc{C}}, \mc{L}, \mc{N}, \varphi, \rho, \mu, \nu) \]
  where
  \begin{enumerate}
  \item $(\mc{C}, \Sigma^{\mc{C}})$ is an $S$-family of genus $g$, $\ell$-pointed, connected twisted curves such that the $i$-th marked point $\Sigma_i$ is banded by $\ab<\gamma_i> \leq \bmu_k$ and has a section $S \to \Sigma_i$;
  \item $\mc{L}$ and $\mc{N}$ are invertible sheaves on $\mc{C}$ such that $(\mc{L},\mc{N})$ is representable, $\deg \mc{L} \otimes \mc{N} = d_0$, $\deg \mc{N} = d_{\infty}$ and the monodromy of $\mc{L}$ along $\Sigma_i^{\mc{C}}$ is $\gamma_i$ when $\ab<\gamma_i> \neq \ab<1>$.
  \item $\mu = (\mu_1, \ldots, \mu_N) \in H^0((\mc{L} \otimes \mc{N})^{\oplus N})$, $\nu \in H^0(\mc{N})$, and $(\mu, \nu)$ is nowhere zero;
  \item $\varphi \in H^0\ab(\bigoplus_{j=1}^N \mc{L}^{\otimes a_j})$ and $(\varphi, \mu)$ is nowhere zero;
    \item $\rho \in H^0(\mc{L}^{-k} \otimes \omega_{\mc{C}/S}^{\log})$ vanishes at $\Sigma_{(1,\rho)}$ and $(\rho, \nu)$ is nowhere zero.
  \end{enumerate}
  We call $\xi$ narrow (resp. broad) if $\gamma$ is narrow (resp. broad).

  An arrow
  \[ \xi' = (\mc{C}'/S', \Sigma^{\mc{C}'}, \mc{L}', \mc{N}', \varphi', \rho', \mu', \nu') \to \xi \]
  consists of a morphism $S' \to S$ and a tuple $(a,b,c)$, where
  \begin{enumerate}
  \item $a \colon (\mc{C}', \Sigma^{\mc{C}'}) \to (\mc{C}, \Sigma^{\mc{C}}) \times_S S'$ is an isomorphism of twisted curves over $S'$;
    \item $b \colon a^* \mc{L} \to \mc{L}'$ and $c \colon a^* \mc{N} \to \mc{N}'$ are isomorphisms of invertible sheaves such that $b^{\otimes a_i}(\varphi_i) = \varphi_i'$, $b^{\otimes -k} \rho = \rho'$, $(b \otimes c)(\mu) = \mu'$, and $c(\nu) = \nu'$.
  \end{enumerate}

  We define $\mc{W}_{g,\gamma,\mbf{d}}^{\pre}$ to be the prestack with the above objects and morphisms.
\end{defn}

\begin{rmk}
If $a_i = 1$ for all $i$, then there the monodromy of $\mc{L}, \mc{N}$ must be opposite at $\Sigma_j^{\mc{C}}$ because either the marked point is a scheme point or $\varphi$ vanishes at $\Sigma_j^{\mc{C}}$, which implies that $\mu$ is nonvanishing at $\Sigma_j^{\mc{C}}$, which trivializes $\mc{L} \otimes \mc{N}$ at $\Sigma_j^{\mc{C}}$. Therefore, representability of $(\mc{L}, \mc{N})$ is equivalent to the representability of both $\mc{L}, \mc{N}$ in this case.

    In our case, for any point $x \in \mc{C}$, either $\mu(x) \neq 0$, which means that $\mc{L}, \mc{N}$ have opposite monodromy at $x$, or $\mu(x) = 0$, in which case $\mc{N}$ has trivial monodromy at $x$ and $\mc{L}|_x$ is a representable line bundle on $\mathbf{B}\!\Aut(x)$. In particular, $\mc{L}$ is always representable and the monodromy of $\mc{N}$ is determined by the type of marking for narrow $\gamma$.
\end{rmk}

\begin{defn}
  Let $S$ be a scheme. We call $\xi \in \mc{W}_{g,\gamma,\mbf{d}}^{\pre}(S)$ \textit{stable} if $\xi|_s$ has finitely many automorphisms for every closed point $s \in S$.
\end{defn}

Having finite automorphism group is an open condition, so we may define $\mc{W}_{g,\gamma,\mbf{d}} \subset \mc{W}_{g,\gamma,\mbf{d}}^{\pre}$ to be the open substack of familes of stable objects. There is an action of $\T = (\C^{\times})^N$ on $\mc{W}_{g,\gamma,\mbf{d}}$ by the formula
\[ t \cdot (\mc{C}, \Sigma^{\mc{C}}, \mc{L}, \mc{N}, \varphi, \rho, \mu, \nu) =(\mc{C}, \Sigma^{\mc{C}}, \mc{L}, \mc{N}, \varphi, \rho, t \cdot \mu, \nu), \]
where for $t = (t_1, \ldots, t_N)$, we define
\[ t \cdot (\mu_1, \ldots, \mu_N) = (t_1 \mu_1, \ldots, t_N \mu_N). \]

\begin{thm}
The stack $\mc{W}_{g,\gamma,\mbf{d}}$ is a Deligne-Mumford $\T$-stack locally of finite type.
\end{thm}

\begin{proof}
  That $\mc{W}_{g,\gamma,\mbf{d}}$ is a Deligne-Mumford stack follows from having finite stabilizers and~\cite[Lemma 2.8]{orbqmap}. To use the cited result, we view $\mc{W}_{g,\gamma,\mbf{d}}$ as the closed substack of $\ul{\Hom}_{\mf{M}_{g,\ell}^{\mr{tw}}}(\mf{C}_{g,\ell}, [\C^{N+7}/(\C^{\times})^3])$ where the third line bundle is $\omega_{\mf{C}_{g,\ell}/\mf{M}_{g,\ell}}^{\log}$. Here, $\mf{C}_{g,\ell}$ denotes the universal curve over $\mf{M}^{\mr{tw}}_{g,\ell}$. That $\mc{W}_{g,\gamma,\mbf{d}}$ is locally of finite type follows from the fact that the stack $\mf{M}_{g,\ell}^{\mr{tw}}$ of prestable twisted curves is locally of finite type~\cite[Corollary 1.12]{logtwcurves}.
\end{proof}

\subsection{Virtual cycle}
\label{subsec:virtual_cycle}

Let $\mc{D}_{g,\gamma}$ be the moduli space of objects $(\mc{C}, \Sigma^{\mc{C}}, \mc{L}, \mc{N})$, where $(\mc{C}, \Sigma_{\mc{C}})$ is a twisted prestable curve and $\mc{L}$ and $\mc{N}$ are invertible sheaves on $\mc{C}$ such that $(\mc{L}, \mc{N})$ is representable and the monodromy of $\mc{L}$ at the marked points is given by $\gamma$. Arrows are given by triples $(a,b,c)$ as in~\Cref{defn:nmsp}. This is a smooth Artin stack of dimension $5g-5+\ell$.

Then let 
\begin{equation}
\pi \colon \mc{C} \to \mc{W}_{g,\gamma,\mbf{d}} 
\end{equation}
be the universal curve over $\mc{W}_{g,\gamma,\mbf{d}}$. Also define 
\begin{equation}
q \colon \mc{W}_{g,\gamma,\mbf{d}} \to \mc{D}_{g,\gamma} 
\end{equation}
to be the forgetful map. For space, we will abbreviate $\mc{P} = \mc{L}^{-k} \otimes \omega_{\mc{C}/\mc{W}_{g,\gamma,\mbf{d}}}^{\log}$.

\begin{prop}
The morphism $q \colon \mc{W}_{g,\gamma,\mbf{d}} \to \mc{D}_{g,\gamma}$ has a relative $\T$-equivariant perfect obstruction theory given by
\begin{equation}
  \ab(R \pi_* \ab(\bigoplus_{i=1}^n \mc{L}^{a_i} \oplus \mc{P}(-\Sigma^{\mc{C}}_{(1,\rho)}) \oplus \bigoplus_{\alpha=1}^N (\mc{L} \otimes \mc{N} \otimes \mbf{L}_{\alpha}) \oplus \mc{N}))^{\vee} \to \L^{\bullet}_{\mc{W}_{g,\gamma,\mbf{d}}/\mc{D}_{g,\gamma}},
\end{equation}
where $\mbf{L}_{\alpha}$ is the character of $\T$ where the $\alpha$-th factor $\T_{\alpha}$ of $\T = (\C^{\times})^n$ acts with weight $1$ and the others act with weight $0$. 
\end{prop}

\begin{proof}
  The proof is the same as in~\cite{mspfermat}.
\end{proof}

The relative obstruction bundle of $q$ is given by
\[ \mc{O}b_{\mc{W}_{g,\gamma,\mbf{d}}/\mc{D}} = R^1 \pi_* \ab(\bigoplus_{i=1}^n \mc{L}^{a_i} \oplus \mc{P}(-\Sigma^{\mc{C}}_{(1,\rho)}) \oplus \bigoplus_{\alpha=1}^N (\mc{L} \otimes \mc{N} \otimes \mbf{L}_{\alpha}) \oplus \mc{N}). \]
The absolute obstruction sheaf is defined by the exact sequence
\[ q^* T_{\mc{D}} \to \mc{O}b_{\mc{W}_{g,\gamma, \mbf{d}}/\mc{D}} \to \mc{O}b_{\mc{W}_{g,\gamma, \mbf{d}}} \to 0. \]
We then define a cosection $\sigma \colon \mc{O}b_{\mc{W}_{g,\gamma,\mbf{d}}/\mc{D}} \to \mc{O}_{\mc{W}_{g,\gamma,\mbf{d}}}$ by the formula
\begin{equation}\label{eqn:cosection}
\sigma(\xi)(\dot{\varphi}, \dot{\rho}, \dot{\mu}, \dot{\nu}) = \sum_{i=1}^n b_i \rho \varphi_i^{b_i - 1} \dot{\varphi}_i + \dot{\rho} \sum_{i=1}^n \varphi_i^{b_i}.
\end{equation}
for any $S$-point $\xi \in \mc{W}_{g,\gamma,\mbf{d}}(S)$.

\begin{lem}
    \Cref{eqn:cosection} defines a $\T$-equivariant homomorphism 
    \[ \sigma \colon \mc{O}b_{\mc{W}_{g,\gamma,\mbf{d}}/\mc{D}} \to \mc{O}_{\mc{W}_{g,\gamma,\mbf{d}}}. \] 
    It lifts to a $\T$-equivariant cosection of $\mc{O}b_{\mc{W}_{g,\gamma,\mbf{d}}}$.
\end{lem}

\begin{proof}
    The proof is the same as in~\cite[Lemma 2.7]{nmsp}.
\end{proof}

Following~\cite{nmsp,mspfermat}, we define the degeneracy locus $\mc{W}^-_{g,\gamma,\mbf{d}}$ with the reduced stack structure. The closed points are characterized by the analogue of~\cite[Lemma 2.13]{mspfermat}.

\begin{lem}
    Let $\gamma$ be narrow. Then $\xi \in \mc{W}_{g,\gamma,\mbf{d}}(\C)$ lies in $\mc{W}_{g,\gamma,\bd}^-(\C)$ if and only if either the locus of $\mc{C}$ where either $\varphi = 0$ or $\sum \varphi_j^{b_j} \rho = 0$ is all of $\mc{C}$.
\end{lem}

\begin{rmk}
    The proceeding discussion also makes sense for $\mc{W}^{\pre}_{g,\gamma,\bd}$ whenever $\gamma$ is narrow, so we define $\mc{W}^{\pre-}_{g,\gamma,\mbf{d}}$ analogously to $\mc{W}^-_{g,\gamma,\mbf{d}}$.
\end{rmk}

Applying the construction of~\cite{cosection}, we obtain a cosection localized virtual cycle
\[ [\mc{W}_{g,\gamma,\mbf{d}}]^{\vir} \in A_*^{\T}(\mc{W}^-_{g,\gamma,\mbf{d}}) \]
with virtual dimension
\[ N d_0 + N(1-g) + d_{\infty} + \ell - 4 \sum_{i=1}^{\ell_0} \frac{m_i}{k}, \]
where $\ell_0$ is the number of stacky insertions and $\gamma_i = \zeta_k^{m_i}$.

\subsection{Definition as $\Omega$-stable moduli space}
\label{subsec:omega stable}

In~\cite{generalmsp}, the authors introduce a generalization of the GLSM moduli spaces introduced in~\cite{glsm}. The basic data is that of
\[ (V, G \leq \Gamma, \epsilon, \vartheta), \]
where $V$ is an affine scheme, $\Gamma$ is a reductive group acting on $V$, $\ep \colon \Gamma \to \C^{\times}$ is a surjective character with kernel $G$, and the stability condition $\vartheta$ is a character of $\Gamma$.

In our case, we set $V = \C^{7+N}$, $G = (\C^{\times})^2 \leq (\C^{\times})^3 = \Gamma$, $\ep(t_1, t_2, t_3) = t_3$, and consider the charge matrix
\[ \begin{pmatrix}
  a_1 & a_2 & a_3 & a_4 & a_5 & -k & 1 & \cdots & 1 & 0 \\
  0 & 0 & 0 & 0 & 0 & 0 & 1 & \cdots & 1 & 1 \\
  0 & 0 & 0 & 0 & 0 & 1 & 0 & \cdots & 0 & 0
\end{pmatrix}. \]
Let $x_1, \ldots, x_5, p, u_1, \ldots, u_N, v$ be the coordinates on $\C^{7+N}$ and define the \textit{master space} by
\[ X = \ab[ \ab(\ab\{(x, u) \neq 0\} \cap \ab\{(u,v) \neq 0\} \cap \ab\{(p,v) \neq 0\}) / G ]. \]
Based on our nonvanishing conditions, we see that the anticones are $(x_i, v)$, $(p,u_j)$, and $(u_j, v)$. Therefore, we need to choose $\vartheta|_G \in \R_{>0} (1,1) + \R_{>0} (0,1)$, for example
\[ \vartheta(t_1, t_2, t_3) = t_1t_2^2. \]

The $\Omega$-stability condition in~\cite{generalmsp} is a triple
\[ \Omega = (S, A, \vartheta), \]
where $S$ is a finite set of elements in
\[ R_+ = \bigoplus_{k > 0} H^0(V, L_{k\vartheta})^G, \]
$A \in \Q$, and $\vartheta$ is as in the basic data.

In our setting, we choose
\[ S = \ab\{x_1 v^{2 a_1}, \ldots, x_5 v^{2 a_5}, uv, u^{2k} p\}\]
and $\frac{1}{k} < A < \frac{2}{k}$.

From this data,~\cite{generalmsp} construct a moduli space
\[ \mr{LGQ}_{g,k}^{\Omega}(X,\mbf{d})\]
of $k$-pointed genus $g$ $\Omega$-stable LG quasimaps. Following~\cite[Appendix A]{generalmsp}, we obtain $\mc{W}_{g,\gamma,\mbf{d}}$ as a connected component of $\mr{LGQ}_{g,k}^{\Omega}(X,\mbf{d})$.

Applying~\cite[Theorem 2.4, Theorem 2.5]{generalmsp}, we obtain
\begin{prop}\label{prop:properness}
    The stack $\mc{W}_{g,\gamma,\mbf{d}}$ is separated and of finite type and the degeneracy locus $\mc{W}^-_{g,\gamma,\mbf{d}}$ is proper.
\end{prop}

\section{$\T$-fixed loci and graphs}
\label{sec:T fixed loci and graphs}

In this section, our goal is to describe the $\T$-fixed loci $\mc{W}_{g,\gamma,\mbf{d}}^{\T}$ in a way similar to~\cite[\S 4]{nmsp} and to~\cite{msp2}. We will fix a narrow $(g,\gamma,\mbf{d})$ and abbreviate $\mc{W}_{g,\gamma,\mbf{d}}$ to $\mc{W}$.

\subsection{Associated graphs}

For each $\xi \in \mc{W}^{\T}(\C)$, we will assign a graph $\Theta_{\xi}$, which has vertices $V$, edges $E$, legs $L$, and flags $F = \ab\{(e,v) \in E \times V \mid v \text{ incident to }e\}$. As always, we will set
\[ \xi = (\mc{C}, \Sigma^{\mc{C}}, \mc{L}, \mc{N}, \varphi,\rho,\mu,\nu). \]
Because $\xi \in \mc{W}^{\T}$, there is an induced action of $\T$ on $\mc{C}$. Let $\mc{C}^{\T}$ be the fixed locus and $\mc{C}^{\mr{mv}} = \mc{C} \setminus \mc{C}^{\T}$. If $\rho$ is nowhere vanishing, we will say $\rho = 1$. The same applies to the other sections.

\begin{defn}
    Let $\xi \in \mc{W}^{\T}$. Define the following unions of components of $\mc{C}$:
    \begin{enumerate}
        \item $\mc{C}_0 = \mc{C} \cap (\mu = 0)_{\mr{red}}$;
        \item $\mc{C}_1 = \mc{C} \cap (\varphi = \rho = 0)_{\mr{red}}$;
        \item $\mc{C}_{\infty} = \mc{C} \cap (\nu = 0)_{\mr{red}}$;
        \item $\mc{C}_{01}$ is the union of irreducible components of $\ol{\mc{C} \setminus (\mc{C}_0 \cup \mc{C}_1 \cup \mc{C}_{\infty})}$ contained in $(\rho = 0)$;
        \item $\mc{C}_{1\infty}$ is the union of irreducible components of $\ol{\mc{C} \setminus (\mc{C}_0 \cup \mc{C}_1 \cup \mc{C}_{\infty})}$ contained in $(\varphi = 0)$;
        \item $\mc{C}_{\infty}$ is the union of irreducible components of $\mc{C}$ not contained in $\mc{C}_0 \cup \mc{C}_1 \cup \mc{C}_{\infty} \cup \mc{C}_{01} \cup \mc{C}_{1\infty}$.
    \end{enumerate}
\end{defn}

If $a,a' \in \Lambda \coloneqq \ab\{0,1,\infty\}$, then $\mc{C}_{aa'}$ consists of curves connecting $\mc{C}_a$ and $\mc{C}_{a'}$.

\begin{defn}
    Let $\xi \in \mc{W}^{\T}$. Then the graph $\Theta_{\xi}$ is defined as follows:
    \begin{enumerate}
        \item\textbf{Vertices } The set $V = V(\Theta_{\xi})$ is the set of connected components of $\mc{C}^{\T}$; 
        \item\textbf{Edges } The set $E$ of edges is the set of connected components of $\mc{C}^{\mr{mv}}$;
        \item\textbf{Legs } The set $L = \ab\{1,\ldots,\ell\}$ is the set of markings, and $i \in L$ is attached to a vertex $v \in V$ if $\Sigma_i^{\mc{C}} \in \mc{C}_v$;
        \item\textbf{Flags } The set $F = \ab\{(e,v) \in E \times V \mid \mc{C}_e \cap \mc{C}_v \neq \emptyset\}$ is the set of flags.
    \end{enumerate}
\end{defn}

There is a partition $V = V_0 \sqcup V_1 \sqcup V_{\infty}$, where $v \in V_a$ if $\mc{C}_v \subseteq \mc{C}_a$ for $a \in \Lambda$. The set of edges can also be partitioned as $E = \bigsqcup_{a,a' \in \Lambda} E_{aa'}$, where $e \in E_{aa'}$ if the two $\T$-fixed points of $\mc{C}_e$ lie in $\mc{C}_a$ and $\mc{C}_{a'}$. From now on, let $\alpha, \beta \in \ab\{1,\ldots,N\}$. We will now further partition the vertices and edges.

\begin{lem}\leavevmode
    \begin{enumerate}
        \item For all $v \in V_{\infty}$, we have $\rho|_{\mc{C}_v} = 1$ and $\nu|_{\mc{C}_v} = 0$. In addition, there exists $\alpha$ such that $\mu_{\alpha}|_{\mc{C}_v} = 0$ and $\mu_{\beta \neq \alpha}|_{\mc{C}_v} = 0$. This induces a partition $V_{\infty} = \bigsqcup_{\alpha} V_{\infty}^{\alpha}$.
        \item For all $v \in V_1$, we have $\varphi|_{\mc{C}_v} = \rho|_{\mc{C}_v} = 0$ and $\nu|_{\mc{C}_v} = 1$. In addition, there exists $\alpha$ such that $\mu_{\alpha}|_{\mc{C}_v} = 0$ and $\mu_{\beta \neq \alpha}|_{\mc{C}_v} = 0$. This induces a partition $V_{1} = \bigsqcup_{\alpha} V_{1}^{\alpha}$.
        \item For all $v \in V_0$, $\varphi|_{\mc{C}_v}$ is nowhere zero, $\mu|_{\mc{C}_v} = 0$, and $\nu|_{\mc{C}_v} = 1$.
    \end{enumerate}
\end{lem}

\begin{proof}
    The statements about $V_0$ follow directly from the definitions of $\mc{C}_0$ and weighted MSP fields.

    If $v \in V$, then $\mc{C}_v \subseteq \mc{C}^{\T}$. If $v \in V_1$, the conditions for $\varphi, \rho, \nu$ follow from the definitions of $\mc{C}_1$ and weighted MSP fields. The statement about $\mu$ follows from ${\T}$-invariance. If $v \in V_{\infty}$, the conditions for $\rho$ and $\nu$ follow from the definitions of $\mc{C}_{\infty}$ and weighted MSP fields, and the statement about $\mu$ follows from $\T$-invariance.
\end{proof}

\begin{lem}
    For all $e \in E_{1\infty} \cup E_{01}$, there exists $\alpha$ such that $\mu_{\alpha}|_{\mc{C}_e \cap \mc{C}_1} \neq 0$ and $\mu_{\beta \neq \alpha}|_{\mc{C}_e \cap \mc{C}_1} = 0$. If $e \in E_{1,\infty}^{\alpha}$, then $\mu_{\alpha}|_{\mc{C}_e} = 1$ and $\mu_{\beta \neq \alpha}|_{\mc{C}_e} = 0$. Therefore, we have partitions $E_{1\infty} = \bigsqcup_{\alpha} E_{1\infty}^{\alpha}$ and $E_{01} = \bigsqcup_{\alpha} E_{01}^{\alpha}$.
\end{lem}

\begin{proof}
    The proof is the same as in~\cite[Lemma 4.4]{nmsp}.
\end{proof}

\begin{conv}
    From now on, we will abbreviate $\mbf{L}_{\alpha} \otimes \mbf{L}_{\beta}^{-1}$ by $\mbf{L}_{\alpha - \beta}$ for all $\alpha, \beta \in \{1, \ldots, N\}$.
\end{conv}

\begin{lem}
    For all $e \in E_{\infty\infty} \cup E_{11}$, there exist $\alpha \neq \beta$ such that $\mu_{\alpha}|_{\mc{C}_e} \neq 0$, $\mu_{\beta}|_{\mc{C}_e} \neq 0$, and $\mu_{\delta \neq \alpha,\beta}|_{\mc{C}_e} = 0$. In addition, $E_{00} = \emptyset$.
\end{lem}

\begin{proof}
    Let $e \in E_{\infty\infty} \cup E_{11}$ and $\ab\{ p, p' \} = \mc{C}_e \cap \mc{C}^{\T}$. We know that $\mu|_{\mc{C}_e}$ is nowhere vanishing, so there exist $\alpha, \beta$ such that $\mu_{\alpha}(p) \neq 0$ and $\mu_{\beta}(p') \neq 0$. If $\alpha = \beta$, then $(\mc{L} \otimes \mc{N})|_{\mc{C}_e} \cong \mc{O}_{\mc{C}_e}$ because $\mu_{\alpha} = 1$. Using the fact that $\deg \mc{L} \leq 0$, we see that all the sections are constant, so $\mc{C}_e$ is unstable. Therefore, $\alpha \neq \beta$. This implies that $\mu_{\alpha}(p') = \mu_{\beta}(p) = 0$ and therefore that $\deg (\mc{L} \otimes \mc{N})|_{\mc{C}_e} > 0$

    We then see that $( \mc{L} \otimes \mc{N} )|_p = \mbf{L}_{\beta - \alpha}$ and $(\mc{L} \otimes \mc{N})|_{p'} = \mbf{L}_{\alpha - \beta}$. This implies that if there exists $\gamma \neq \alpha,\beta$ such that $\mu_{\gamma}|_{\mc{C}_e} \neq 0$, then $\mu_{\gamma}(p) = \mu_{\gamma}(p') = 0$. But then if $\T_{\gamma}$ acts nontrivially on $\mc{C}_e$, this would imply that $\deg ( \mc{L} \otimes \mc{N} )|_{\mc{C}_e} = 0$, which is a contradiction. This implies that $\T_{\gamma}$ acts trivially on $\mc{C}_e$, so $\mu_{\gamma} = 0$.

    The statement that $E_{00} = \emptyset$ follows from the fact that $\mc{C}_0 \subseteq \mc{C}^{\T}$, which follows from the fact that $\mu|_{\mc{C}_0} = 0$.
\end{proof}

\begin{defn}\label{defn:decorations}
    We will add decorations to $\Theta_{\xi}$ as follows:
    \begin{enumerate}
        \item \textbf{Genus } For $a \in V \cup E$, we set $g_a = H^1(\mc{C}_a, \mc{O}_{\mc{C}_a})$;
        \item \textbf{Degree } For $a \in V \cup E$, we set $\mbf{d}_a = (\deg(\mc{L} \otimes \mc{N})|_{\mc{C}_a}, \deg \mc{N}|_{\mc{C}_a}) \eqqcolon (d_{0a}, d_{\infty a})$;
        \item \textbf{Monodromy } The monodromy of $\mc{L}$ at a leg $i \in \ab\{1,\ldots,\ell\}$ is that of the monodromy of the marking it represents (recall that $\mc{L}$ is always representable). Because we are in the narrow sector, the monodromy of $\mc{N}$ is uniquely determined by the type of marking.
        \item \textbf{Monodromy } Monodromy is assigned to flags as follows:
            \begin{enumerate}
                \item If $e \in E_{01}$ and $v \in V_0$, then if $\deg \mc{L}|_{\mc{C}_e} = a + \frac{1}{k}b$ for $b \in [1,r]$, we assign $\gamma_{(e,v)} = (\zeta_k^b, \rho)$. Note if $b \neq r$, then the degeneracy locus is empty, so it suffices to consider graphs where we assign $(1,\rho)$;
                \item If $e \in E_{01}$ and $v \in V_1$, then we assign $\gamma_{(e,v)} = (1, \rho)$;
                \item If $e \in E_{1\infty}$ and $v \in V_1$, then we assign $\gamma_{(e,v)} = (1, \varphi)$;
                \item If $e \in E_{1\infty}$ and $v \in V_{\infty}$, then if $\deg \mc{L}|_{\mc{C}_e} = a + \frac{1}{k}b$ for $b \in [1,k]$, we assign $\gamma_{(e,v)} = (\zeta_k^b, \varphi)$;
                \item If $e \in E_{0\infty}$ and $v \in V_0$, then if $\deg (\mc{L} \otimes \mc{N})|_{\mc{C}_e} = a + \frac{1}{k}b$ for $b \in [1,k]$, we assign $\gamma_{(e,v)} = \zeta_k^b$;
                \item If $e \in E_{0\infty}$ and $v \in V_0$, then if $\deg \mc{N}|_{\mc{C}_e} = a + \frac{1}{k}b$ for $b \in [1,k]$, we assign $\gamma_{(e,v)} = \zeta_k^b$. Again, if $b \neq r$, the degeneracy locus is empty, so it suffices to consider graphs where we assign $1$.
            \end{enumerate}
            Again, the monodromy of $\mc{N}$ is uniquely determined by the level of the vertex.
        \item \textbf{Hour } Elements of $V_1^{\alpha}$ and $V_{\infty}^{\alpha}$ are assigned hour $\alpha$.
        \item \textbf{Level } Elements of $V_{\infty}$ and $E_{\infty\infty}$ have level $\infty$, elements of $V_1$ and $E_{11}$ have level $1$, and elements of $V_0$ have level $0$.
    \end{enumerate}
\end{defn}

We adopt the following conventions. We call a vertex $v \in V$ \textit{stable} if $\dim \mc{C}_v = 1$ and \textit{unstable} otherwise and denote the set of stable vertices by $V^S \subset V$. For any $v \in V$, we let $E_v$ (resp. $L_v$) be the set of edges (resp. legs) incident to $v$. For any $v \in V^S$ and $e \in E_v$, we call the associated node $q_{(e,v)} \coloneqq \mc{C}_e \cap \mc{C}_v$. Finally, we also decompose
\[ V \setminus V^S = V^{0,1} \cup V^{1,1} \cup V^{0,2}, \]
where $V^{a,b}$ is the set of vertices with $a$ legs and $b$ edges incident to them.

\begin{warn}
    There may exist $\xi_1$ and $\xi_2$ in the same connected component of $\mc{W}^{\T}$ where $\Theta_{\xi_1} \not\cong \Theta_{\xi_2}$.
\end{warn}

\subsection{Flattening}

Let $G_{g,\gamma,\mbf{d}}$ be the set of all isomorphism classes of decorated graphs $\Theta_{\xi}$ of $\xi \in \mc{W}^{\T}(\C)$. We will discuss how to modify $\Theta_{\xi}$ in a way to produce a unique ``flattened'' graph $\Theta_{\xi}^{\mr{fl}}$ describing each connected component of $\mc{W}^{\T}$.

\begin{defn}
    Let $\mc{C}$ be a $\T$-twisted curve and $q$ be a node of $\mc{C}$. Then the formal completion of $\mc{C}$ at $q$ is a union $\wh{\mc{C}}_q = \wh{\mc{C}}_1 \cup \wh{\mc{C}}_2$ of two irreducible components. Then we define $q$ to be \textit{$\T$-balanced} if $T_q \wh{\mc{C}}_1 \otimes T_q \wh{\mc{C}}_2 \cong \mbf{L}_0$.
\end{defn}

All nodes $q_{(e,v)}$ of $\xi$ coming from flags $(e,v) \in F$ where $v \in V^S$ are $\T$-unbalanced, with $T_{q_{(e,v)}} \wh{\mc{C}}_e$ nontrivial and $T_{q_{(e,v)}} \wh{\mc{C}}_v$ trivial. The others are given by the set $V^{0,2}$.
These $v$ are unstable with $\ab|L_v| = 0$ and $\ab|E_v| = 2$. We denote their associated nodes by $q_v$. 

\begin{lem}\label{lem:balancednode}
    Let $v \in V^{0,2}$ and let $e, e'$ be the two edges incident to $v$. Let $d_e = \deg \mc{L}|_{\mc{C}_e}$ and $d_{e'} = \deg \mc{L}|_{\mc{C}_{e'}}$ Then $q_v$ is $\T$-balanced if and only if either $(\mc{C}_e \cup \mc{C}_{e'}) \cap \mc{C}_{\infty}$ is a special point and $d_e + d_{e'} = 0$.
\end{lem}

\begin{proof}
    Suppose that $v \in V_{\infty}$. If $e \in E_{\infty\infty}$ and $e' \in E_{1\infty}$, then if $\alpha$ is the hour of $v$, $\T_{\alpha}$ is the only factor of $\T$ acting nontrivially on $\mc{C}_{e'}$. However, there exists $\beta \neq \alpha$ such that $\T_{\beta}$ acts nontrivially on $\mc{C}_e$, so $v$ is not balanced. The next case is when $e' \in E_{\infty\infty}$. Then there exists a unique $\alpha$ such that $\mu_{\alpha}|_{q_v} \neq 0$. Then either there exists $\beta$ such that $e,e' \in E_{\alpha\beta}$, in which case $T_{q_v} \mc{C}_e$ and $T_{q_v} \mc{C}_{e'}$ are positive multiples of each other, or $e$ and $e'$ are in different $E_{\alpha\beta}$, in which case $q_v$ is clearly not balanced. If $e' \in E_{0\infty}$, there exists a unique $\gamma$ such that $\mu_{\gamma}(q_v) \neq 0$ by $\T$-invariance and the fact that $\nu$ is nonzero on $\mc{C}_{e'}$. But then the diagonal subgroup of $\T$ acts nontrivially on $\mc{C}_{e'}$, so $q_v$ cannot be balanced.

    The next case is when $e \in E_{1\infty}$. If $e' \in E_{1\infty}$, they must have the same hour. Then $T_{q_v} \mc{C}_e$ and $T_{q_v} \mc{C}_{e'}$ are positive multiples of each other, so $q_v$ is not balanced. If $e' \in E_{0\infty}$, then by the same argument as above, we see that $q_v$ is not balanced.

    Finally, if $e, e' \in E_{0\infty}$, then $T_{q_v} \mc{C}_e$ and $T_{q_v} \mc{C}_{e'}$ must be positive multiples of each other.

    If $v \in V_0$, then the cases where $e, e' \in E_{01}$, where $e \in E_{01}$ and $e' \in E_{0\infty}$, and where $e, e' \in E_{0\infty}$ are treated exactly as above.

    Finally, let $v \in V_1$ and let $\alpha$ be the hour of $v$. The cases where $e \in E_{11}$, $e,e' \in E_{01}$ or $e,e' \in E_{1\infty}$ are treated as above. Therefore, we assume that $e \in E_{01}$ and $e' \in E_{1\infty}$. Let $q = \mc{C}_e \cap \mc{C}_0$ and $q' = \mc{C}_{e'} \cap \mc{C}_{\infty}$. Because $\varphi(q) \neq 0$ and $\rho'(q') \neq 0$, we have 
    \[ \mc{L}|_q \cong (\mc{L}^{-k} \otimes \omega_{\mc{C}}^{\log})|_{q'} \cong \mbf{L}_0. \] 
    Also, because $q_v$ is a node, we know $\mc{L}^{-k}|_{q_v} \cong (\mc{L}^{-k} \otimes \omega^{\log}_{\mc{C}})|_{q_v}$ as $\T_{\alpha}$-representations (recall that $\T_{\alpha}$ is the only factor of $\T$ acting nontrivially on $C_e$ and $C_{e'}$). Studying the $\T_{\alpha}$-action on $\mc{C}_e$ and $\mc{C}_{e'}$, we observe that $q_v$ is balanced if and only if
    \[ \deg \mc{L}^{-k} |_{\mc{C}_e} + \deg (\mc{L}^{-k} \otimes \omega^{\log}_{\mc{C}})|_{\mc{C}_{e'}} = 0. \]
    This is the same as 
    \[ -k d_e - k d_{e'} + \deg \omega^{\log}_{\mc{C}}|_{\mc{C}_{e'}} = 0. \]
    The desired result follows immediately from the observation that we cannot have $\deg \omega^{\log}_{\mc{C}}|_{\mc{C}_{e'}} = -1$ because otherwise both $\mc{C}_e$ and $\mc{C}_{e'}$ are schemes, so we cannot have $k d_e + k d_{e'} = -1$.
\end{proof}

\begin{defn}
    A graph $\Theta \in G_{g,\gamma,\mbf{d}}$ is \textit{flat} if it has no $\T$-balanced nodes.
\end{defn}

\begin{defn}
    Let $\Theta \in G_{g,\gamma,\mbf{d}}$. Then we construct the \textit{flattening} $\Theta^{\mr{fl}}$ of $\Theta$ by the following procedure: for every $\T$-balanced $v \in N(\Theta)$ with edges $e \in E_{01}$ and $e' \in E_{1\infty}$ incident to it, we eliminate $v$ from $\Theta$ and combine $e, e'$ into a single edge $\wt{e} \in E_{0\infty}$. The decorations are $g_{\wt{e}} = 0$, $d_{0\wt{e}} = d_{0e}$, and $d_{\infty \wt{e}} = d_{\infty e'}$ (note that $(\mc{L} \otimes \mc{N})|_{\mc{C}_{e'}} \cong \mc{O}_{\mc{C}_{e'}}$ and $\mc{N}|_{\mc{C}_e} \cong \mc{O}_{\mc{C}_e}$). Monodromy is assigned via~\Cref{defn:decorations}.
\end{defn}

Let $G^{\fl}_{g,\gamma,\mbf{d}} = \ab\{\Theta^{\fl} \mid \Theta \in G_{g,\gamma,\mbf{d}}\}/\sim$. For any $\Theta \in G^{\fl}_{g,\gamma,\mbf{d}}$, we define a $\Theta$-framed MSP field to be a pair $(\xi, \epsilon)$, where $\epsilon \colon \Theta_{\xi}^{\fl} \cong \Theta$ is an isomorphism of decorated graphs. Following~\cite[\S 4.2]{nmsp}, we can form the stack of $\Theta$-framed MSP fields with the morphism
\[ \iota_{\Theta} \colon \mc{W}_{\Theta} \to \mc{W}^{\T}. \]
Denote its image by $\mc{W}_{(\Theta)}$. Because $\mc{W}_{\Theta} \to \mc{W}_{(\Theta)}$ is an $\Aut \Theta$-torsor, we have the equality
\[ [W_{\Theta}]^{\vir} = \ab|\Aut \Theta| \cdot [\mc{W}_{(\Theta)}]^{\vir}. \]
In addition, there is a decomposition
\[ \mc{W}^{\T} = \bigsqcup_{\Theta \in G_{g,\gamma,\mbf{d}}^{\mr{fl}}} \mc{W}_{(\Theta)}. \]

\section{Irregular vanishing}
\label{sec:irrvan}

We will now prove that the virtual cycles for some graphs vanish. This will greatly simplify our study of the virtual localization formula.

\begin{defn}
    A graph $\Theta$ is a \textit{pure loop} if it has no stable vertices and every vertex has exactly two edges attached to it.
\end{defn}

\begin{defn}
  A vertex $v \in V_{\infty}$ is \textit{regular} if either
  \begin{enumerate}
    \item If $v$ is stable, then $\gamma_v = (1^{m_1}2^{m_2})$ where $2$ appears only if $\frac{2}{k}$ is narrow;
    \item If $v$ is unstable attached to an edge in $E_{1\infty}$, then $\mc{C}_v$ is a nonscheme marking.
  \end{enumerate}
  A flat graph $\Theta$ is \textit{regular} if all of its vertices $v \in V_{\infty}$ are regular and $E_{0\infty} = \emptyset$. If it is not regular and not a pure loop, we call it \textit{irregular}.
\end{defn}

\begin{rmk}
    The condition on the monodromy of the insertions at level $\infty$ is determined by numerical criteria for the nonvanishing of FJRW invariants, see~\cite{fjrw} for details.
\end{rmk}

\begin{defn}
    Given a cycle class $A \in A_*^{\T} \mc{W}^-$, we say that $A \sim 0$ if there is a proper $\T$-invariant substack $\mc{W}^- \subset \mc{W}' \subset \mc{W}$ such that $A = 0$ when viewed as elements of $A_*^{\T}(\mc{W}')$.
\end{defn}

\begin{thm}\label{thm:irrvan}
  If $\Theta \in G_{g,\gamma,\mbf{d}}^{\mr{fl}}$ is irregular, then $[\mc{W}_{\Theta}]^{\vir} \sim 0$.
\end{thm}

We will now prove this theorem following the proof in~\cite[\S 5]{nmsp}. Because our case is very similar to the case of the quintic, we will content ourselves with an outline of the proof. In this section, if we cite a result from~\cite{msp3} or~\cite{nmsp} in the statement of a lemma or theorem, it means that the proof is the same as the proof of the cited result. 

The key points to check are that
\begin{enumerate}
  \item In the ``no-string'' case, the virtual dimension is negative;
  \item If $\bar{e}$ is a single $0\infty$ edge (a ``string''), then $\mc{W}_{\bar{e}}$ has at worse local complete intersection singularities.
\end{enumerate}
When these two facts are satisfied, the argument in~\cite{nmsp} also works in our case.

\subsection{Alternative constructions of $[ \mc{W}_{\Theta} ]^{\vir}$}%
\label{sub:Alternative constructions of virtual cycle}

We begin with some alternative constructions of the MSP virtual cycle for the fixed loci. 
\begin{defn}
    A \textit{$\Theta$-framed twisted curve} is a $\T$-equivariant $(\mc{C}, \Sigma^{\mc{C}}, \mc{L}, \mc{N}) \in \mc{D}$ together with an identification $\ep$ of the following data that is consistent with the geometry of $\mc{C}$:
    \begin{itemize}
        \item The marked points $\Sigma^{\mc{C}}$ are identified with the legs of $\Theta$;
        \item The $\T$-unbalanced nodes of $\mc{C}$ are identified with the nodes of $\Theta$;
        \item Let $\mc{C}^{\mr{dec}}$ be the normalization of $\mc{C}$ along its $\T$-unbalanced nodes. We will identify the connected components of $\mc{C}^{\mr{dec}}$ with $V^S(\Theta) \cup E(\Theta)$.
    \end{itemize}
    We require that when $e \in E_{0\infty}(\Theta)$, then either $\mc{C}_e \cong \P^1$ or $\mc{C}_e = \mc{C}_{e-} \cup \mc{C}_{e+}$ is a union of two copies of $\P^1$ such that $\mc{C}_{e-} \cap \mc{C}_0 \neq \emptyset$, $\mc{C}_{e+} \cap \mc{C}_{\infty} \neq \emptyset$, $\mc{L}\otimes \mc{N} \otimes \mbf{L}_{\alpha}|_{\mc{C}_{e+}} \cong \mc{O}_{\mc{C}_{e+}}$ for some $\alpha = 1, \ldots, N$, and $\mc{N}|_{\mc{C}_{e-}} \cong \mc{O}_{\mc{C}_{e-}}$.
\end{defn}

We will denote the stack of $\Theta$-framed twisted curves by $\mc{D}_{\Theta}$. This is a smooth Artin stack.

\begin{defn}
    A \textit{$\Theta$-framed gauged twisted curve} is an object $(\mc{C}, \Sigma^{\mc{C}}, \mc{L}, \mc{N}, \mu, \nu)$ which is $\T$-equivariant together with an identitification $\ep$ such that
    \begin{enumerate}
        \item $(\mc{C}, \Sigma^{\mc{C}}, \mc{L}, \mc{N}, \ep) \in \mc{D}_{\Theta}$ is a $\Theta$-framed twisted curve;
        \item $\mu \in H^0\ab(\bigoplus_{\alpha=1}^N \mc{L}\otimes \mc{N} \otimes \mbf{L}_{\alpha})^{\T}$ and $\nu \in H^0(\mc{N})^{\T}$ such that
            \begin{enumerate}
                \item $\mu|_{\mc{C}_0} = 0$, $\nu|_{\mc{C}_{\infty}} = 0$, and $\mu|_{\mc{C}_{\infty}}$ and $\nu|_{\mc{C}_0}$ are nowhere vanishing;
                \item For every $e \in E_{1\infty}$, there exists $\alpha \in \{1, \ldots, N \}$ such that  $\mu_{\alpha}|_{\mc{C}_e}$ is nowhere vanishing and $\mu_{\beta\neq\alpha}|_{\mc{C}_e} = 0$;
                \item For every $e \in E_{01}$, there exists $\alpha \in \{1, \ldots, N \}$ such that  $\mu_{\alpha}|_{\mc{C}_e}\neq 0$ and $\mu_{\beta\neq\alpha}|_{\mc{C}_e} = 0$;
                \item If $e \in E_{0\infty}$ is such that $\mc{C}_e$ is a union of two copies of $\P^1$, then $\nu|_{\mc{C}_{e-}}$ is nowhere vanishing and the condition for $\mu$ in (b) holds for $\mc{C}_{e+}$.
            \end{enumerate}
    \end{enumerate}
\end{defn}

The stack $\mc{D}_{\Theta, \nu}$ of all $\Theta$-framed gauged twisted curves is a smooth Artin stack. As in~\cite[\S 3]{msp3}, we have morphisms
\[ \mc{W}_{\Theta} \to \mc{D}_{\Theta, \nu} \to \mc{D}_{\Theta} \]
and relative perfect obstruction theories for $\mc{W}_{\Theta}$ over $\mc{D}_{\Theta, \nu}$ and $\mc{D}_{\Theta}$. Denote the corresponding virtual fundamental classes by $[\mc{W}_{\Theta}]^{\vir}_{\Theta, \nu}$ and $[\mc{W}_{\Theta}]^{\vir}_{\Theta}$.

\begin{lem}[{\cite[Proposition 3.4]{msp3}}]
    We have an equality
    \[ [\mc{W}_{\Theta}]^{\vir} = [\mc{W}_{\Theta}]^{\vir}_{\Theta, \nu} = [\mc{W}_{\Theta}]^{\vir}_{\Theta}. \]
\end{lem}

\subsection{Trimming the graph}%
\label{sub:Trimming edges}

We will now simplify the analysis by removing edges from $\Theta$ in a way that preserves the vanishing of the virtual cycle. We will then study the remaining graphs first in the case where $\Theta$ has no ``string'' by reducing the $N=1$ case, and then when $\Theta$ has a string we will use a virtual pullback argument relying on the geometric structure of $\mc{W}_{\bar{e}}$.

We will need the notion of \textit{webs}, which were first introduced for the quintic in~\cite[\S 5.2]{nmsp}. For any $\Theta \in G_{g,\gamma,\mbf{d}}^{\mr{fl}}$, let $\Theta_{\infty}$ be the (possibly disconnected) graph obtained by deleting all vertices in $V_0 \cup V_1$ (and all legs attached to them) and edges in $E_{01} \cup E_{11}$ and replacing all edges in $E_{1\infty} \cup E_{0\infty}$ with legs. These legs are assigned monodromy as in~\Cref{defn:decorations}. Some components of $\Theta_{\infty}$ are degenerate, having no stable vertices and no edges.

\begin{defn}
    A connected nondegenerate component of $\Theta_{\infty}$ is called a \textit{web} of $\Theta$.
\end{defn}

We will now explain how to obtain a graph $\Theta'$ such that $E_{01}(\Theta') = E_{1\infty}(\Theta') = V_1(\Theta') = \emptyset$.
\begin{enumerate}
    \item We will simply remove all edges in $E_{01}(\Theta)$ and $E_{1\infty}(\Theta)$ and all vertices in $V_1(\Theta)$;
    \item Next, we discard all degenerate connected components. These come from unstable vertices of $\Theta$;
    \item Finally, when we remove an edge in $E_{1\infty}$, we replace it with a leg as in the construction of webs. When we remove an edge in $E_{01}$, we replace it with a $(1,\rho)$ marked point.
\end{enumerate}

\begin{lem}
    In the scenario above, if $[\mc{W}_{\Theta'}]^{\vir} \sim 0$, then
    \[ [\mc{W}_{\Theta}]^{\vir} \sim 0. \]
\end{lem}

\begin{proof}
    This follows from the argument used to prove~\cite[Lemma 3.11]{msp2}.
\end{proof}

Even in this case, we need to make further simplifications by trimming $\frac{1}{k}$ legs and some edges and vertices. Note that in the FJRW theory CohFT, the $\frac{1}{k}$ insertion is the unit and therefore satisfies the string equation~\cite[Theorem 5.1.2]{mfcohftFJRW}.

\begin{defn}
    Let $\Theta$ satisfty $E_{01}(\Theta) = E_{1\infty}(\Theta) = V_1(\Theta) = \emptyset$.
    \begin{enumerate}
        \item Let $e \in E_{\infty\infty}(\Theta)$ be incident to vertices $v_-$ and $v_+$. Suppose that $q_{(e,v_+)}$ is a node in $\mc{C}$. Then 
            \begin{itemize}
                \item $e$ is of \textit{marking type} if $v_- \in V^{1,1}$;
                \item $e$ is of \textit{leaf-end type} if $v_- \in V^{0,1}$;
                \item $e$ is of \textit{nodal type} if $q_{(e,v_-)}$ is a node in $\mc{C}$;
            \end{itemize}
        \item A leg $\ell \in L(\Theta)$ is \textit{spare} if it is decorated by $\frac{1}{k}$ and is attached to a vertex $v \in V(\Theta)$ such that $2g_v + n_v > 3$;
        \item A vertex $v \in V(\Theta)$ is \textit{neutral} if both $E_v$ and $L_v$ are nonempty, $\on{val}(v) = 3$, and at least one leg attached to $v$ is a $\frac{1}{k}$ marking.
    \end{enumerate}
\end{defn}

We will trim all edges of marking or leaf-end type. Let $\mf{a}$ be a web of $\Theta$. 
\begin{itemize}
    \item If $e \in E(\mf{a})$ is of marking type let $\ell$ be the leg attached to $v_-$. We remove $e$ and $v_-$ from $\mf{a}$ and attach $\ell$ to $v_+$ with unchanged decoration;
    \item If $e$ is of leaf-end type, because $\deg \omega_{\mc{C}_e}^{\log} = -1$, the fact that $\rho|_{\mc{C}_e}$ is nowhere vanishing implies that $\deg \mc{L}|_{\mc{C}_e} = -\frac{1}{k}$. This implies that the node of $\mc{C}_e$ has monodromy $\frac{k-1}{k}$. We will delete $e$ from $\mf{a}$ and attach a leg decorated by $\frac{1}{k}$ to $v_+$.
\end{itemize}
We will denote the result of this process by $\mf{a}_e$.

When $\ell$ is a spare marking, we may remove $\ell$ from $\Theta$ and the resulting graph $\Theta'$ is still a graph of $\T$-fixed MSP fields.

Finally, we consider the case when $v \in V(\Theta)$ is a neutral vertex. In this case, let $\ell$ be the leg with the $\frac{1}{k}$ marking. Then because $g_v = 0$ and $\rho|_{\mc{C}_v}$ is nowhere vanishing, we must have $\deg \mc{L}|_{\mc{C}_v} = \frac{1}{k}$, which implies that the monodromy of $\mc{L}$ at the other two special points are equal to each other. We construct $\Theta'$ as follows:
\begin{itemize}
    \item If $v$ has one edge and two legs, let $\ell'$ be decorated by $\frac{a}{k}$. We remove $v$ from $\Theta$ and replace the edge $e$ by a new leg decorated by $\frac{a}{k}$;
    \item If $v$ has one leg and two edges, then we delete the leg and make $v$ unstable.
\end{itemize}
\begin{lem}
    In either case, there is a canonical morphism $\psi \colon \mc{W}_{\Theta} \to \mc{W}_{\Theta'}$ such that
    \[ [\mc{W}_{\Theta}]^{\vir} =  \psi^* [\mc{W}_{\Theta'}]^{\vir}. \]
\end{lem}

\begin{proof}
    The proof follows from the same arguments as in the proofs of~\cite[Lemmas 5.7--5.10]{nmsp}.
\end{proof}

\subsection{The no-string case}
\label{subsec:no string}

\begin{defn}
    Let $\Theta \in G_{g,\gamma,\bd}^{\mr{fl}}$. A \textit{string} of $\Theta$ is an edge $e \in E_{0\infty}(\Theta)$ such that the vertex at level $0$ is unstable and has no other edge attached to it.
\end{defn}

\begin{lem}
    Let $\Theta \in G_{g,\gamma,\bd}^{\mr{fl}}$. Suppose $\Theta$ does not contain any strings. Then
    \[ [\mc{W}_{\Theta}]^{\vir} \sim 0. \]
\end{lem}

\begin{proof}
    By the previous discussion, we may assume that $\Theta$ has no vertices at level $1$, edges of marking or leaf-end type, spare legs, or neutral vertices. Any $\frac{1}{k}$ marking must be attached to an unstable vertex $v$ that has an edge $e \in E_{0\infty}$ attached to it by assumption. But then we may assume that $\mc{W}_{\Theta}^-$ is nonempty, so $e$ must arise from flattening a balanced node. But~\Cref{lem:balancednode} implies that $\mc{C}_e \cap \mc{C}_{\infty}$ is a scheme point, which is a contradiction.

    The next step is to prove that the virtual dimension is negative in this case. The main tool is to reduce to the case when $N = 1$, which requires us to reduce all webs $\mf{a}$ into one-vertex graphs and then delete all hour decorations from $\Theta$. We will let $\mf{a}'$ be a one-vertex graph whose vertex $v$ whose decorations are as follows:
    \begin{itemize}
        \item The genus $g$ is the total genus of $\mf{a}$;
        \item The legs and edges incident to $\mf{a}'$ are the legs and edges incident to $\mf{a}$;
        \item All decorations of legs and edges remain the same, except we delete the hours.
        \item Note that every one-vertex graph of MSP fields with $N=1$ at level $\infty$ must satisfy $\mc{L} \otimes \mc{N} \cong \mc{O}_{\mc{C}}$, so $d_0 = 0$. Using this and the fact that $\rho$ is everywhere nonzero at level $\infty$ (implying that $\mc{L}^{-k} \otimes \omega_{\mc{C}}^{\log} \cong \mc{O}_{\mc{C}}$), there is a unique choice of $d_{\infty}$ making $\mf{a}'$ a graph of MSP fields with $N=1$.
    \end{itemize}

    We may apply the above procedure to all webs of $\Theta$ and create a graph $\Theta_1$ of MSP fields with $N=1$. By the argument in the proof of~\cite[Lemma 5.12]{nmsp}, we see that
    \[ \vir.\dim \mc{W}_{\Theta} \leq \vir.\dim \mc{W}_{\Theta_1}. \]
    Therefore, we may assume that $N=1$.

    In this case, we follow the arguments in~\cite[\S 4]{msp3}. All of the virtual dimension computations are the same as in the quintic case because we are in the Calabi-Yau setting, so we see that
    \[ \vir.\dim \mc{W}_{\Theta_1} < 0, \]
    which implies the desired result.
\end{proof}

\subsection{The general case}%
\label{sub:general irrvan}

We will prove~\Cref{thm:irrvan} by reduction to the no-string case. Let $e$ be a string of $\Theta$. Then we will define the graph $\bar{e}$ as follows the graph with one edge $e$ and two vertices $v_- \in V_0$ and $v_+ \in V_{\infty}$ together with the decorations on $e$, a leg decorated by $1 = \zeta_k^0$ at $v_+$, and possibly a leg at $v_-$. Then let $\mc{W}_{\bar{e}}$ be the moduli space of stable $\bar{e}$-framed MSP fields.

\begin{rmk}
    Note that the stability condition for MSP fields implies that any $0\infty$ edge must come from flattening a balanced node. Indeed, any geometric $0\infty$ edge must satisfy $\deg \mc{N} > 0$ and $\deg \mc{L}^{-k} \otimes \omega_{\mc{C}_e}^{\log} > 0$, which implies that $\deg \mc{L} < 0$ and hence $\varphi \equiv 0$. But this is impossible because $\mu = 0$ at level $0$ and the stability condition requires $(\varphi, \mu) \neq 0$ everywhere.
\end{rmk}

\begin{lem}
    The stack $\mc{W}_{\bar{e}}$ has pure dimension $4$, hypersurface singularities, and virtual dimension $4$.
\end{lem}

\begin{proof}
    We will give a construction of $\mc{W}_{\bar{e}}$. For simplicity, we will assume that $N=1$, but the general case is the same, except we introduce hours. Let $R = \Bl_{\infty \times 0} \P^1 \times \A^1$ and $E_{\infty} \subset R$ be the exceptional divisor. Also let $E_0$ be the strict transform of $\P^1 \times 0$, $D_0$ be the strict transform of $0 \times \A^1$, and $D_{\infty}$ be the strict transform of $\infty \times \A^1$. Let $u$ be the coordinate on $\A^1$ and $[x,y]$ be the homogeneous coordinates on $\P^1$. Finally, set $\delta = d_{0e}$.

    Now let $X = (u^{k\delta} = 0) \subset \A^1$ and set $\mc{C}_X \coloneqq R \times_{\A^1} X$ and $\Sigma^{\mc{C}_X} = D_{\infty}$. We now introduce sections. Set
    \[ L = \mc{O}_R(\delta E_{\infty}) \qquad \text{and} \qquad  N = \mc{O}_R(\delta D_{\infty}). \]
    Then we define the sections
    \begin{align*}
        \varphi &\colon \mc{O}_R \to L = \mc{O}_R(\delta E_{\infty}) \\
        \mu &\colon \mc{O}_R \to L \otimes N = \mc{O}_R(\delta E_{\infty} + \delta E_{\infty}) = \mc{O}_R(\delta D_0) \\
        \nu &\colon \mc{O}_R \to N = \mc{O}_R(\delta D_{\infty})
    \end{align*}
    to be the defining inclusions and denote their restrictions to $\mc{C}_X$ by $\varphi_X$, $\mu_X$, and $\nu_X$, respectively.

    We need to be more care about $\rho$. Here, we set
    \[ P \coloneqq \omega_{R/\A^1}(D_{\infty}) = \mc{O}_R(-D_0 - k\delta E_{\infty}). \]
    This has a section $\rho_0$ on the central fiber given by the following commutative diagram.
    \begin{equation*}
    \begin{tikzcd}
        \mc{O}_R(-D_0)|_{R_0} \ar{r}{u^{k\delta}} \ar{d} & \mc{O}_R(-D_0 - k\delta E_0 - k \delta)|_{R_0} \ar[hookrightarrow]{r} & \mc{O}_R(-D_0 - k\delta E_{\infty})|_{R_0} \ar[equal]{d} \\
        \mc{O}_{R_0} \ar{rr}{\rho_0} & & P|_{R_0}.
    \end{tikzcd}
    \end{equation*}
    Note that $\rho$ cannot extend to all of $R$ because any section of $P$ must vanish outside the central fiber, but it does extend to $X$. In fact, $X$ is the maximal subscheme $S \subset \A^1$ such that $\rho$ extends to $S$. We denote the resulting lift by $\rho_X$.

    We now give $(\mc{C}_X, \Sigma^{\mc{C}_X}, L|_{\mc{C}_X}, N|_{\mc{C}_X}, \varphi_X, \rho_X, \mu_X, \nu_X)$ a $\T = \C^{\times}$ action. First, we let $\T$ act on $\P^1 \times \A^1$ by
    \[ t \cdot ([x,y], u) \coloneqq ([tx,y],u). \]
    This action lifts to an action on $R$ which then restricts to $\mc{C}_X$. We then linearize $\mc{O}_R$ such that its restriction to $\mc{O}_{D_0}$ is the trivial action and give $L$, $N$, and $P$ the $\T$-linearizations which make the defining inclusions $\mc{O}_R \subset L$, $\mc{O}_R \subset N$, and $P \subset \mc{O}_R$ $\T$-equivariant.

    We then note that we must have
    \[ (\varphi_X, \rho_X, \mu_X, \nu_X)^t = (t^{b_1} \varphi_X, t^{b_2} \rho_X, t^{b_3} \mu_X, t^{b_4} \nu_X) \]
    for some rational number $b_j$. We know that $\varphi_X$ and $\nu_X$ are nowhere vanishing along $D_0$, so we must have $b_1 = b_4 = 0$. Next, $\mu_X$ has vanishing order $\delta$ along $D_0$, so $b_2 = \delta$. The definition of $\rho_X$ and equivariance of $P \subset \mc{O}_R$ imply that $b_2 = 0$.

    The upshot of the previous discussion is that $\mc{W}_{\bar{e}}$ is not reduced. We now construct an explicit model for it. Let $\mc{X} \coloneqq (\A^5 \setminus 0) \times X$ and $\mc{C}_{\mc{X}} \coloneqq (\A^5 \setminus 0) \times \mc{C}_X$. We will give this family a section
    \[ \Sigma^{\mc{C}_{\mc{X}}} \coloneqq (\A^5 \setminus 0) \times D_{\infty} \times_X \mc{C}_X \]
    with the monodromy $1 = \zeta_k^0$. Let $x_1, \ldots, x_5$ be the coordinates on $\A^5$ and $q \colon \mc{C}_{\mc{X}} \to \mc{C}_X$ be the natural projection. We then define the line bundles
    \[ \mc{L} \coloneqq q^* L \qquad \text{and}  \qquad\mc{N} \coloneqq q^* N \]
    and the sections
    \[ \varphi_i \coloneqq x_i \cdot q^* \varphi^{a_i}_X, \qquad \rho = q^* \rho_X, \qquad \mu = q^* \mu_X, \qquad \text{and} \qquad \nu = q^* \nu_X. \]
    This is a family of $\bar{e}$-framed MSP fields, so we obtain a morphism
    \[ \eta \colon \mc{X} \to \mc{W}_{\bar{e}}. \]

    By the proof of~\cite[Lemma 8.1]{msp3}, $\eta$ factors through the $\C^{\times}_{\tau}$-action on $\mc{X}$ given by
    \[ \tau \cdot (x_1, x_2, x_3, x_4, x_5, u) \coloneqq (\tau^{a_1\delta} x_1, \ldots, \tau^{a_5 \delta}, \tau^{-1}u) \]
    and induces an isomorphism $\mf{X} / \C^{\times} \cong \mc{W}_{\bar{e}}$. Because $X \subset \A^1$ was defined by a single equation, $\mc{W}_{\bar{e}}$ has hypersurface singularities.

    We finally prove the statement about the virtual dimension. We note that $H^1(\mc{L})^{\T} = 0$ and
    \[ \dim H^0(\mc{L}^{a_i})^{\T} = \dim H^0(\mc{P})^{\T} =\dim H^1(\mc{P})^{\T} = 1, \]
    so the relative dimension of $\mc{W}_{\bar{e}}$ over $\mc{D}_{\bar{e}, \nu}$ is
    \[ \sum_{i=1}^5\chi_{\T}(\mc{L}^{a_i}) + \chi_{\T}(\mc{P}) = 5. \]
    Beceause $\dim \mc{D}_{\bar{e}, \nu} = -1$ by the formulae in~\cite[\S 4]{msp3}, the virtual dimension of $\mc{W}_{\bar{e}}$ is $4$.
\end{proof}

We may now follow the arguments in~\cite{msp3} to conclude the proof of~\Cref{thm:irrvan}.

\section{Virtual localization}
\label{sec:Virtual localization}

We will state the virtual localization formula. Even though the proof is the same as in~\cite{msp2,nmsp}, we will give some of the arguments to clarify the exposition in~\cite{nmsp} and correct a slight mistake in~\cite[Lemma 4.4]{msp2}. Let $\Theta$ be a regular graph with webs $\mf{a}_1, \ldots, \mf{a}_s$. Also, let $Z \subset \P(\ba)$ be the Fermat Calabi-Yau hypersurface.

\subsection{Fixed parts}%
\label{sub:Fixed parts}

Note that
\[ W_{\Theta} \cong \ab( \prod_{v \in V_0} \mc{W}_v \times_{\P(\ba)^{\abs{E_v}}} \prod_{e \in E_v}\mc{W}_e ) \times \prod_{v \in V_1^S} \mc{W}_v \times \prod_{k=1}^s \mc{W}_{\mf{a}_k} \times \prod_{e \in E_{1\infty} \cup E_{11}} \mc{W}_e. \]
Now define $G_e$ to be the automorphism group of a point in $\mc{W}_e$ for any $e \in E$ and set
\[ G_E \coloneqq \prod_{e \in E \setminus E_{\infty\infty}} G_e. \]
In the diagram
\begin{equation*}
\begin{tikzcd}
    \mc{W}_e \ar{r}{\iota_{\Theta}} \ar{d} & \mc{W}_{(\Theta)} \subset \mc{W}^{\T} \\
    \prod_{v \in V_0 \cup V_1^S} \mc{W}_v \times \prod_{k=1}^s \mc{W}_{\mf{a}_k},
\end{tikzcd}
\end{equation*}
note the top arrow is a $\Aut(G)$-torsor and the vertical arrow is a $G_e$-gerbe. Using~\cite[Theorem 3.5]{locwccosection}, we have
\begin{prop}
    Let $\Theta \in G_{g,\gamma,\mbf{d}}^{\mr{reg}}$. Then the virtual class of the fixed component corresponding to $\Theta$ is
    \[ [W_{(\Theta)}]^{\vir} = \frac{1}{\ab|\Aut(\Theta)|} \frac{1}{\ab|G_E|} (\iota_{\Theta})_* \prod_{v \in V_0 \cup V_1^S} [\mc{W}_v]^{\vir} \times \prod_{m=1}^s [\mc{W}_{\mf{a}_m}]^{\vir}. \]
\end{prop}

\begin{lem}
    The cosection localized virtual cycles are as follows:
    \begin{enumerate}
        \item For $v \in V_0^S$, $\mc{W}_v$ is the GLSM moduli space for the geometric phase and we have
            \[ [\mc{W}_v]^{\vir} = (-1)^{1-g_v} [\ol{\mc{M}}_{g,E_v \cup \Sigma_v}(Z, d_v)]^{\vir}; \]
        \item For $v \in V_1^S$, $\mc{M}_v = \ol{\mc{M}}_{g_v, E_v \cup \Sigma_v}$ and $[\mc{W}_v]^{\vir} = [\mc{W}_v]$.
        \item For a web $\mf{a}$ at infinity, all legs of $\mf{a}$ are narrow by irregular vanishing, and $[\mc{W}_{\mf{a}}]^{\vir} = [\mc{W}_{(\Theta')}]^{\vir}$ for $\Theta' = \mf{a}$.
        \item For $v \in V_0^U$, then $[\mc{W}_v]^{\vir} = -[Z]$.
    \end{enumerate}
\end{lem}

\subsection{Moving parts}%
\label{sub:Moving parts}

We now write down the moving parts. While the proof is the same as in~\cite{msp2,nmsp}, we will give slightly more detail than in~\cite[\S 6]{nmsp} and correct some slight mistakes in the formulae relating to the presence of unstable vertices.

Let
\[ \mc{V} \coloneqq \bigoplus_{i=1}^5 \mc{L}^{a_i} \oplus \mc{P}(-\Sigma^{\mc{C}}_{(1,\rho)}) \oplus \bigoplus_{\alpha=1}^N (\mc{L} \otimes \mc{N} \otimes \mbf{L}_{\alpha}) \oplus \mc{N}. \]
Following~\cite[\S 4]{msp2}, denote
\begin{align*}
    B_1 &= \Aut(\mc{C}, \Sigma^{\mc{C}}) = \Ext^0(\Omega^{\mc{C}}(\Sigma^{\mc{C}}), \mc{O}_{\mc{C}}), \\
    B_2 &= \Aut(\mc{L}) \oplus \Aut(\mc{N}) = H^0(\mc{O}_{\mc{C}}^{\oplus 2}), \\
    B_3 &= \on{Def}(\varphi, \rho, \mu, \nu) = H^0(\mc{V}) \\
    B_4 &= \on{Def}(\mc{C}, \Sigma^{\mc{C}}) = \Ext^1(\Omega_{\mc{C}}(\Sigma_{\mc{C}}), \mc{O}_{\mc{C}}) \\
    B_5 &= \on{Def}(\mc{L}) \oplus \on{Def}(\mc{N}) = H^1(\mc{O}_{\mc{C}}^{\oplus 2}) \\
    B_6 &= \on{Obs}(\varphi, \rho, \mu, \nu) = H^1(\mc{V}).
\end{align*}
Using this notation, we have
\[ \frac{1}{e_{\T}(N^{\vir})} = \frac{e_{\T}(B_1^{\mv}) e_{\T}(B_2^{\mv}) e_{\T}(B_6^{\mv})}{e_{\T}(B_3^{\mv}) e_{\T}(B_4^{\mv}) e_{\T}(B_5^{\mv})}. \]

In addition, recall the normalization exact sequences
\begin{align*}
    0 \to \mc{V} \to \bigoplus_{v \in V^S \cup E} \mc{V}|_{\mc{C}_v} \to \bigoplus_{a \in F^S \cup V^{0,2}} \mc{V}|_a \to 0, \\
    0 \to \mc{O}_{\mc{C}} \to \bigoplus_{v \in V^S \cup E} \mc{O}_{\mc{C}_v} \to \bigoplus_{a \in F^S \cup V^{0,2}} \mc{O}_{\mc{C}}|_a \to 0. 
\end{align*}

\begin{notn}
    For $\alpha, \beta = 1, \ldots, N$, denote $\mbf{L}_{\alpha} \otimes \mbf{L}_{\beta}^{-1}$ by $\mbf{L}_{\alpha-\beta}$ and $\mbf{L}_{\alpha}^{k}$ by $\mbf{L}_{k\alpha}$ for any $k \in \Q$.
\end{notn}

\begin{notn}
    For $a \in V^S \cup E \cup F^S \cup V^{0,2}$, we will let $A'_a$ denote the contribution of $a$ to $\frac{1}{e_{\T}(N^{\vir})}$.
\end{notn}

\subsubsection{Contribution from stable vertices}%
\label{ssub:Contribution from stable vertices}

For $v \in V^S$, let $\pi_v \colon \mc{C}_v \to \mc{W}_v$ be the universal curve and $\mc{L}_v$ be the universal line bundle corresponding to $\mc{L}$. Also, let $\E_v = (\pi_v)_* \omega_{\pi_v}$ be the Hodge bundle.

First, suppose $v \in V_0$. Then let  $f_v \colon \mc{C} \to \P(\ba)$ be the universal morphism defined by $(\mc{L}, \varphi)$. Note that $\nu|_{\mc{C}_v} = 1$ (implying $\mc{N}|_{\mc{C}_v} \cong \mc{O}_{\mc{C}_v}$) and deformations of $\varphi, \rho$ are in the fixed part, so the moving part is
\[ (H^0(\mc{V}|_{\mc{C}_v}) - H^1(\mc{V}|_{\mc{C}_v}))^{\mv} - (H^0(\mc{O}_{\mc{C}_v}^{\oplus 2}) - H^1(\mc{O}_{\mc{C}_v}^{\oplus 2}))^{\mv} = H^0\ab(\bigoplus_{\alpha=1}^N \mc{L} \otimes \mc{N} \otimes \mbf{L}_{\alpha} | _{\mc{C}_v}). \]
Therefore, we obtain
\[ A_v' = \frac{1}{\prod_{\alpha=1}^N e_T(R(\pi_v)_* f_v^* \mc{O}_{\P(\ba)}(1) \otimes \mbf{L}_{\alpha})}. \]

Now suppose $v \in V_1^{\alpha}$. Then $\varphi = \rho = 0$, $\mu_{\alpha} = 1$, $\mu_{\beta \neq \alpha} = 0$, and $\nu = 1$, so the moving part is
\[ (\chi(\mc{V}|_{\mc{C}_v}))^{\mv} - (\chi(\mc{O}_{\mc{C}_v}^{\oplus 2}))^{\mv} = \chi \ab(\bigoplus_{i=1}^5 \mc{L}^{a_i}|_{\mc{C}_v} \oplus \mc{P}(-\Sigma_{(1,\rho)}^{\mc{C}})|_{\mc{C}_v} \bigoplus_{\beta \neq \alpha}^N \mc{L} \otimes \mc{N} \otimes \mbf{L}_{\beta}|_{\mc{C}_v}). \]
Using the fact that 
\[ \mc{L} \otimes \mc{N} \otimes \mbf{L}_{\alpha} |_{\mc{C}_v} \cong \mc{N}|_{\mc{C}_v} \cong \mc{O}_{\mc{C}_v}, \]
we see that
\begin{align*}
    \mc{L}|_{\mc{C}_v} &= \mbf{L}_{-\alpha}, \\
    \mc{P}(-\Sigma_{(1,\rho)}^{\mc{C}})|_{\mc{C}_v} &= \mbf{L}_{r \alpha} \otimes \omega_{\mc{C}}|_{\mc{C}_v}, \\
    \mc{L}\otimes \mc{N} \otimes \mbf{L}_{\beta} &= \mbf{L}_{\beta - \alpha}.
\end{align*}
Taking the Euler characteristics and then the Euler classes, we obtain
\[ A_v' = \prod_{i=1}^5 \frac{e_T(\E_v^{\vee} \otimes \mbf{L}_{-a_i \alpha})}{-a_i t_{\alpha}} \cdot \frac{kt_{\alpha}}{e_T(\E_v \otimes \mbf{L}_{k\alpha}) \cdot (k t_{\alpha})^{\abs{E_v}}} \cdot \prod_{\beta \neq \alpha} \frac{e_T(\E_v^{\vee} \otimes \mbf{L}_{\beta - \alpha})}{t_{\beta} - t_{\alpha}}. \]

Now let $v \in V_{\infty}^{\alpha}$. Note that $\rho = 1$, $\mu_{\alpha} = 1$, $\mu_{\beta \neq \alpha} = 0$, and $\nu = 0$. Also, deformations of $\varphi$ are fixed, so the moving part is
\[ \chi(\mc{V}|_{\mc{C}_v})^{\mv} - \chi(\mc{O}_{\mc{C}_v})^{\mv} = \chi\ab(\bigoplus_{\beta \neq \alpha} \mc{L}\otimes \mc{N} \otimes \mbf{L}_{\beta} |_{\mc{C}_v} \oplus \mc{N} |_{\mc{C}_v}). \]
Because
\begin{align*}
    \mc{N}|_{\mc{C}_v} &= \mc{L}^{-1}|_{\mc{C}_v} \otimes \mbf{L}_{-\alpha}, \\
    \mc{L}\otimes \mc{N} \otimes \mbf{L}_{\beta} &= \mbf{L}_{\beta - \alpha},
\end{align*}
we obtain
\[ A_v' = \frac{1}{e_T(R (\pi_v)_* \mc{L}_v^{-1} \otimes \mbf{L}_{-\alpha})} \prod_{\beta \neq \alpha} \frac{e_T(\E_v^{\vee} \otimes \mbf{L}_{\beta - \alpha})}{t_{\beta} - t_{\alpha}}. \]

\subsubsection{Edge contributions}%
\label{ssub:Edge contributions}

For $e \in E_{01} \cup E_{1\infty}$, let $v, v'$ be the incident vertices such that $v' \in V_1^{\alpha}$ and let $q, q' \in \mc{C}_e$ be the corresponding points. Define
\begin{align*}
    \delta = \begin{cases}
        -1 & v \in V^{0,1} \\
        0 & \text{otherwise},
    \end{cases} \qquad
    \delta' = \begin{cases}
        -1 & v' \in V^{0,1} \\
        0 & \text{otherwise}.
    \end{cases}
\end{align*}

First, let $e \in E_{01}$. Note that $\nu = 1$ while $\mu_{\alpha}(q) = \varphi(q') = 0$ and $\mu_{\alpha}(q'), \varphi(q) \neq 0$. Also, note that $d_e > 0$. Let
\[ \mc{E}_e \coloneqq \mc{O}_{\P(d_e \cdot \ba)}(1) \to \mc{W}_e = \P(d_e \cdot \ba) \]
be the universal bundle and $h_e$ be the hyperplane class on $\P(\ba)$. Also define
\begin{align*}
    \delta_{\rho} = \begin{cases}
        -1 & v \in V^{0,1} \\
        0 & \text{otherwise},
    \end{cases} \qquad
    \delta_{\rho}' = \begin{cases}
        -1 & v \in V^{0,1} \\
        0 & \text{otherwise}.
    \end{cases}
\end{align*}
As $\T_{\alpha}$-equivariant line bundles, we have
\begin{align*}
    \mc{L}|_{\mc{C}_e} &= \mc{O}_{\mc{C}_e}(d_e q') \\
    \mc{P}(-\Sigma_{(1,\rho)}^{\mc{C}}) &= \mc{O}_{\mc{C}_e}((\delta + \delta_{\rho}) q + (-k d_e + \delta' + \delta_{\rho}')q') \\
    \mc{L} \otimes \mc{N} \otimes \mc{L}_{\alpha} &= \mc{O}_{\mc{C}_e}(d_e q).
\end{align*}
Standard computational techniques (for example, see~\cite[Example 98]{localization}) give us
\begin{align*}
    R^0 (\pi_e)_* \mc{L}^{a_i}|_{\mc{C}_e} &= \bigoplus_{j=0}^{a_i d_e} \mc{E}_e^{a_i d_e - j} \otimes \mbf{L}_{-\frac{j}{d_e} \alpha} \\
    R^1 (\pi_e)_* \mc{P}(-\Sigma_{(1,\rho)}^{\mc{C}})|_{\mc{C}_e} &= \bigoplus_{j=1+\delta+\delta_{\rho}}^{k d_e - 1 - \delta - \delta'_{\rho}} \mc{E}_e^{-kd_e + \delta' + \delta'_{\rho} + j} \otimes \mbf{L}_{\frac{j}{d_e} \alpha} \\
    R^0 (\pi_e)_* \mc{L} \otimes \mc{N} \otimes \mbf{L}_{\alpha} |_{\mc{C}_e} &= \bigoplus_{j=0}^{d_e} \mc{E}_e^j \otimes \mbf{L}_{\frac{j}{d_e} \alpha} \\
    R^0 (\pi_e)_* \mc{L} \otimes \mc{N} \otimes \mbf{L}_{\beta} |_{\mc{C}_e} &= \bigoplus_{j=0}^{d_e} \mc{E}_e^j \otimes \mbf{L}_{\frac{j}{d_e} \alpha} \otimes \mbf{L}_{\beta - \alpha}. 
\end{align*}
Taking equivariant Euler classes, we obtain
\[ A'_e = \frac{\prod_{j=1}^{kd_e-1-\delta'-\delta'_{\rho}} \ab(-k h_e +\frac{\delta'+\delta'_{\rho}}{d_e} h_e + j \frac{h_e + t_{\alpha}}{d_e})}{\prod_{i=1}^5 \prod_{j=1}^{a_i d_e} \ab(a_i h_e - j\frac{h_e + t_{\alpha}}{d_e}) \prod_{j=1}^{d_e} j\frac{h_e+t_{\alpha}}{d_e} \prod_{\beta \neq \alpha} \prod_{j=1}^{d_e} \ab(j\frac{h_e + t_{\alpha}}{d_e} + t_{\beta} - t_{\alpha})}. \]

Now let $e \in E_{1\infty}^{\alpha}$ and set $k_e \coloneqq \max\ab\{ r \in \Z>0 \mid d_e \in \frac{1}{k}\Z \}$. First, suppose $v \notin V^{0,1}$. Note that $d_e < 0$, $\varphi = 0$, $\mu_{\alpha} = 1$, $\mu_{\beta \neq \alpha} = 0$, and the zeroes of $\rho$ are concentrated at $q'$ while the zeroes of $\nu$ are concentrated at $q$. This implies that
\begin{align*}
    \mc{L} \otimes \mc{N} \otimes \mc{L}_{\alpha}|_{\mc{C}_e} &= \mc{O}_{\mc{C}_e} \\
    \mc{N}|_{\mc{C}_e} &= \mc{O}_{\mc{C}_e}(k_e d_e q) \\
    \mc{P}(-\Sigma_{(1,\rho)}^{\mc{C}}) \ab((-kd_e + \delta' + \delta'_{\rho})q')
\end{align*}
as $\T_{\alpha}$-equivariant line bundles, so we obtain
\begin{align*}
    H^1(\mc{L}^{a_i}|_{\mc{C}_e}) &= \bigoplus_{j=1}^{\ceil{-a_i d_e}-1} \mbf{L}_{-\ab(a_i + \frac{j}{d_e})\alpha} \\
    H^0(\mc{N}|_{\mc{C}_e}) &= \bigoplus_{j=0}^{\floor{- d_e}} \mbf{L}_{\frac{j}{d_e} \alpha} \\
    H^0(\mc{P}(-\Sigma_{(1,\rho)}^{\mc{C}})|_{\mc{C}_e}) &= \bigoplus_{j=k d_e - \delta' - \delta'_{\rho}}^0 \mbf{L}_{\frac{j}{d_e}\alpha} \\
    H^0(\mc{L}\otimes \mc{N}\otimes \mbf{L}_{\beta}) &= \mbf{L}_{\beta-\alpha}.
\end{align*}
We therefore obtain
\[ A'_e = \frac{\prod_{i=1}^5 \prod_{j=1}^{\ceil{-a_i d_e}-1} \ab(-a_i t_{\alpha} + \frac{j}{d_e}t_{\alpha})}{\prod_{j=1}^{-k d_e + \delta' + \delta'_{\rho}} \frac{-j}{d_e} t_{\alpha} \prod_{j=1}^{\floor{-d_e}} \frac{j}{d_e}t_{\alpha} \prod_{\beta \neq \alpha} (t_{\beta}-t_{\alpha})}. \]

Now suppose $v \in V^{0,1}$. Then $\mc{C}_e$ is a scheme and $\delta = -1$. By the argument in the proof of~\cite[Lemma 4.1]{msp2}, we see that the tangent weights $w_{(e,v)}, w_{(e,v')}$ at $q, q'$, respectively are given by
\[ w_{(e,v)} = \frac{k t_{\alpha}}{k d_e+1}, \qquad w_{(e,v')} = \frac{-k t_{\alpha}}{k d_e + 1}. \]
By the same calculations as above, we obtain
\[ A'_e = \frac{\prod_{i=1}^5 \prod_{j=1}^{-a_i d_e-1} \ab(-a_i t_{\alpha} + \frac{kj}{kd_e+1}t_{\alpha})}{\prod_{j=1}^{-k d_e -1 + \delta' + \delta'_{\rho}} \frac{-r j}{k d_e+1} t_{\alpha} \prod_{j=1}^{\floor{-d_e}} \frac{kj}{kd_e+1}t_{\alpha} \prod_{\beta \neq \alpha} (t_{\beta}-t_{\alpha})}. \]
Unifying the two cases, we obtain
\[ A'_e = \frac{\prod_{i=1}^5 \prod_{j=1}^{\ceil{-a_i d_e}-1} \ab(-a_i t_{\alpha} + \frac{kj}{kd_e - \delta}t_{\alpha})}{\prod_{j=1}^{-k d_e + \delta' + \delta'_{\rho}} \frac{-kj}{kd_e-\delta} t_{\alpha} \prod_{j=1}^{\floor{-d_e}} \frac{kj}{kd_e - \delta}t_{\alpha} \prod_{\beta \neq \alpha} (t_{\beta}-t_{\alpha})} \]
for any $e \in E_{1\infty}^{\alpha}$.

The final case is when $e \in E_{11}^{\alpha\beta}$. Recall that $v' \in V_1^{\alpha}$ and $v \in V_1^{\beta}$. Note that $\mu_{\alpha}(v) = \mu_{\beta}(v') = 0$, while $\nu = 1$ is fixed. As $\T_{\alpha-\beta}$-equivariant line bundles, we have
\begin{align*}
    \mc{L} \otimes \mc{N} \otimes \mbf{L}_{\alpha} |_{\mc{C}_e} &= \mc{O}_{\mc{C}_e}(d_e q) \\
    \mc{L} \otimes \mc{N} \otimes \mbf{L}_{\beta} |_{\mc{C}_e} &= \mc{O}_{\mc{C}_e}(d_e q') \\
    \mc{L} |_{\mc{C}_e} &= \mc{O}_{\mc{C}_e}(d_e q) \otimes \mbf{L}_{-\alpha} \\
    \mc{P}(-\Sigma_{(1,\rho)}^{\mc{C}}) |_{\mc{C}_e} &= \mc{O}_{\mc{C}_e}(( -k d_e +\delta + \delta_{\rho} )q + (\delta' + \delta'_{\rho})q') \otimes \mbf{L}_{k\alpha},
\end{align*}
so taking the relevant cohomology groups gives us
\begin{align*}
    H^0(\mc{L}\otimes\mc{N}\otimes\mbf{L}_{\alpha}|_{\mc{C}_e}) &= \bigoplus_{j=0}^{d_e} \mbf{L}_{j\frac{\alpha-\beta}{d_e}} \\
    H^0(\mc{L}\otimes\mc{N}\otimes\mbf{L}_{\beta}|_{\mc{C}_e}) &= \bigoplus_{j=0}^{d_e} \mbf{L}_{j\frac{\beta-\alpha}{d_e}} \\
    H^0(\mc{L}\otimes\mc{N}\otimes\mbf{L}_{\gamma}|_{\mc{C}_e}) &= \bigoplus_{j=0}^{d_e} \mbf{L}_{j\frac{\alpha-\beta}{d_e} + \gamma-\alpha} \\
    H^0(\mc{L}^{a_i}|_{\mc{C}_e}) &= \bigoplus_{j=0}^{a_i d_e} \mbf{L}_{j\frac{\alpha-\beta}{d_e} -a_i \alpha} \\
    H^1(\mc{P}(-\Sigma_{(1,\rho)}^{\mc{C}})|_{\mc{C}_e}) &= \bigoplus_{j=1+\delta'+\delta'_{\rho}}^{k d_e - 1 - \delta - \delta_{\rho}} \mbf{L}_{k \alpha -  \delta -\delta_{\rho}  - j\frac{\alpha-\beta}{d_e}},
\end{align*}
so taking the Euler classes of the moving parts yields
\[ A_e' =  \frac{(-1)^{d_e}(d_e)^{2 d_e}}{(d_e!)^2 (t_{\beta}-t_{\alpha})^{2d_e}} \frac{\displaystyle\prod_{j=1+\delta'+\delta'_{\rho}}^{k d_e - 1 - \delta - \delta_{\rho}}\ab(k t_{\alpha} - \frac{\delta+\delta_{\rho}}{d_e} t_{\alpha} - j \frac{t_{\alpha}-t_{\beta}}{d_e})}{\displaystyle\prod_{\substack{i=1 \\ a+b = a_i d_e}}^5  \ab(-\frac{a}{d_e} t_{\alpha} - \frac{b}{d_e} t_{\beta}) \prod_{\substack{ \gamma \neq \alpha,\beta \\ a+b = d_e }}  \ab(t_{\gamma} - \frac{a}{d_e} t_{\alpha} - \frac{b}{d_e} t_{\beta})}. \]

\subsubsection{Node contributions}%
\label{ssub:Node contributions}

We now study the contributions of nodes in $F^S \cup V^{0,2}$. First, let $(e,v) \in F^S$. The first case is when $v \in V_0$. Then $\mc{W}_v$ is the moduli space for the geometric phase of the GLSM given by $V = \C^6$ and $G = \C^{\times}$ acting with weights $a_1, \ldots, a_5, -r$ and $R$-charges $0,\ldots,0,1$ (for the definition of a GLSM, see~\cite{glsm}). Then let $\on{ev}_{(e,v)}$ be the evaluation map at the marking labeled by $e$. We see that $\varphi$, $\rho$, and $\nu=1$ are fixed, so the moving part is
\[ (H^0(\mc{V}|_{q_{(e,v)}}) - H^0(\mc{O}_{q_{(e,v)}}^{\oplus 2}))^{\mv} = \on{ev}_{(e,v)}^* \mc{O}_{\P(\ba)}(1) \otimes \mbf{L}_{\alpha}. \]
When $v \in V_1^{\alpha}$, $\mu_{\alpha} = \nu = 1$ are fixed and everything else is moving. Then by the same calculations as before, the moving part is
\[ (H^0(\mc{V}|_{q_{(e,v)}}) - H^0(\mc{O}_{q_{(e,v)}}^{\oplus 2}))^{\mv} = \bigoplus_{i=1}^5 \mbf{L}_{-a_i \alpha} \oplus \mbf{L}_{r \alpha} \oplus \bigoplus_{\beta \neq \alpha} \mbf{L}_{\beta-\alpha}. \]
Finally, let $v \in V_{\infty}$ be stacky and narrow (in the non-narrow case, the formula is more complicated). Then $\rho = \mu_{\alpha}=1$ while $\mc{L}^{a_i}, \mc{N}$ have nontrivial monodromy at $q_{(e,v)}$, so 
\[ H^0(\mc{V}|_{q_{(e,v)}}) - H^0(\mc{O}_{q_{(e,v)}}^{\oplus 2}) = \bigoplus_{\beta \neq \alpha} \mbf{L}_{\beta - \alpha}. \]
Therefore, we have
\[ A'_{(e,v)} = \begin{cases}
    \prod_{\alpha=1}^N (h_e + t_{\alpha}) & v \in V_0 \\
    - k \prod_{i=1}^5 a_i t_{\alpha}^6 \prod_{\beta \neq \alpha} (t_{\alpha} - t_{\beta}) & v \in V_1^{\alpha} \\
    \prod_{\beta \neq \alpha} (t_{\beta} - t_{\alpha}) & v \in V_{\infty}^{\alpha} \text{ is narrow in FJRW theory}.
\end{cases}
\]

Under the same assumptions (that if $v \in V_{\infty}^{0,2}$, then it is narrow in FJRW theory), we obtain the same answers for $v \in V^{0,2}$ by the argument used to prove~\cite[Lemma 4.8]{msp2}.

\subsubsection{Contributions from deforming $(\mc{C}, \Sigma^{\mc{C}})$}%

We now compute the remaining contributions to $\frac{1}{e_T(N^{\vir})}$, namely
\[ \frac{e_T(B_1^{\mv})}{e_T(B_4^{\mv})}. \]
We follow the argument in~\cite[\S 4.1]{msp2}. Let $e$ be an edge and $v, v'$ be the two vertices of an edge. As before, by convention we set $v \in V_1^{\alpha}$. Define
\[ w_{(e,v)} \coloneqq \T_{q_{(e,v)}} \mc{C}_v. \]
Formulae are given by the following:
\begin{enumerate}
    \item If $v \in V_0$, then $w_{(e,v)} = \frac{h_e + t_{\alpha}}{d_e}$ and $w_{(e,v')} = -\frac{h_e+t_{\alpha}}{d_e}$;
    \item If $v \in V_{\infty}^{\alpha} \setminus V^{0,1}$, then $w_{(e,v)} = \frac{t_{\alpha}}{k_e d_e}$ and $w_{(e,v')} = -\frac{t_{\alpha}}{d_e}$;
    \item If $v \in V_{\infty}^{\alpha} \cap V^{0,1}$, then $w_{(e,v)} = \frac{kt_{\alpha}}{k d_e+1}$ and $w_{(e,v')} = -\frac{kt_{\alpha}}{k d_e+1}$;
    \item If $e \in E_{11}^{\alpha\beta}$, then $v \in V_1^{\beta}$. In this case, $w_{(e,v)} = \frac{t_{\alpha}-t_{\beta}}{d_e}$.
\end{enumerate}
Then the localization contributions are
\begin{align*}
    e_T(B_1^{\mv}) &= \prod_{(e,v) \in F^{(0,1)}} w_{(e,v)} \\
    e_T(B_4^{\mv}) &= \prod_{(e,v) \in F^S} (w_{(e,v)} - \psi_{(e,v)}) \prod_{\substack{ v \in V^{0,2} \\ E_v = \ab\{e,e'\} }} (w_{(e,v)} + w_{(e',v)}).
\end{align*}

\subsubsection{Contribution from webs at level $\infty$}%
\label{ssub:Contribution from webs at level infty}

Let $\mf{a}_1, \ldots, \mf{a}_s$ be the webs of $\Theta$. 
Set
\[ A_{\infty} = \frac{1}{\prod_{i=1}^s e_T(N_{\mf{a}_i}^{\vir})}. \]
We will see later that this is the remaining contribution after we account for all terms outside of level $\infty$.

\subsection{The virtual localization formula}%
\label{sub:The virtual localization formula}

Define the miscellaneous contributions from the remaining unstable vertices as follows. First, for $v \in V$, define
\[ A_v = \begin{cases}
    A'_v & v \in V^S \cup V^{1,1} \\
    \frac{A'_v}{w_{(e,v)} + w_{(e',v)}} & v \in V^{0,2} \\
    w_{(e,v)} & v \in V^{0,1}.
\end{cases}
\]
Also, for a flag $(e,v) \in F^S$, define
\[ A_{(e,v)} = \frac{A'_{(e,v)}}{w_{(e,v)} - \psi_{(e,v)}}. \]
In the previous subsection, we proved the following result.

\begin{prop}
    The Euler class of the virtual normal bundle $e_T(N_{\Theta}^{\vir})$ is given by
    \[ \frac{1}{e_T(N_{\Theta}^{\vir})} = \ab( \prod_{\substack{(e,v) \in F^S \\ e \notin E_{\infty\infty}}} A_{(e,v)} ) \cdot \ab( \prod_{v \in V \setminus V_{\infty}} A_v ) \cdot \ab( \prod_{e \in E \setminus E_{\infty\infty}} A'_e ) \cdot A_{\infty}. \]
\end{prop}

\appendix

\section{Alternative proof of properness}
\label{sec:properties}

In this appendix, we give an alternative proof of~\Cref{prop:properness} in line with that in~\cite{mspfermat,nmsp}. In the process, we fix some details in~\cite{mspfermat}.

\subsection{Stability criterion for the degeneracy locus}
\label{subsec:stability}

\begin{lem}\label{lem:stability}
    Let $\xi \in \mc{W}^{\pre}$. Let 
    \[\pi \colon \wt{\mc{C}} = \bigsqcup_{\alpha} \wt{\mc{C}}_\alpha \to \mc{C} \] 
    be the normalization of $\mc{C}$, where the $\wt{\mc{C}}_\alpha$ are the connected components of $\wt{\mc{C}}$. Let $\Sigma^{\wt{\mc{C}}} = \pi^{-1}(\Sigma^{\mc{C}} \cup \mc{C}_{\mr{sing}})$, and let $(\wt{\mc{C}}, \wt{\mc{N}}, \wt{\varphi}, \wt{\rho}, \wt{\mu}, \wt{\nu})$ be the pullbacks of $(\mc{L}, \mc{N}, \varphi, \rho,\mu,\nu)$ via $\pi^*$. Setting
    \[ \wt{\xi}_\alpha = (\wt{\mc{C}}_\alpha, \Sigma^{\wt{\mc{C}}} \cap \wt{\mc{C}}_\alpha, \wt{\mc{L}}|_{\wt{\mc{C}}_\alpha},\wt{\mc{N}}|_{\wt{\mc{C}}_\alpha},\wt{\varphi}|_{\wt{\mc{C}}_\alpha},\wt{\rho}|_{\wt{\mc{C}}_\alpha},\wt{\mu}|_{\wt{\mc{C}}_\alpha},\wt{\nu}|_{\wt{\mc{C}}_\alpha}), \]
    then $\xi$ is stable if and only if all $\wt{\xi}_\alpha$ are stable.
\end{lem}

\begin{proof}
  By the same argument as in~\cite[Lemma 3.3]{mspfermat}, if we define $\Aut_{\mc{E}}(\xi)$ to be the subgroup of $\Aut(\xi)$ that fixes $\mc{E}$ for any component $\mc{E}$ of $\mc{C}$, $\mc{C}$ must contain an irreducible component $\mc{C}_\alpha$ such that the image of $\Aut_{\mc{C}_\alpha}(\xi) \to \Aut(\mc{C}_\alpha)$ is infinite. But this image is exactly $\Aut(\wt{\xi}_\alpha)$.
\end{proof}

\begin{lem}\label{lem:stabilitydegeneracy}
  Let $\xi \in \mc{W}^{\pre -}(\C)$. Then $\xi$ is unstable if and only if one of the following holds:
  \begin{enumerate}
  \item $\mc{C}$ contains a rational curve $\mc{E}$ which contains one special point and either $\rho|_{\mc{E}}$ is nowhere zero and $\mc{L} \otimes \mc{N} |_{\mc{E}} \cong \mc{O}_{\mc{E}}$ or $\mc{L}^k|_{\mc{E}} \cong \mc{O}_{\mc{E}}$ and $\mc{N}|_{\mc{E}} \cong \mc{O}_{\mc{E}}$;
  \item $\mc{C}$ contains a rational curve $\mc{E}$ which contains two special points, $\mc{L}^k |_{\mc{E}} \cong \mc{O}_{\mc{E}}$, and either $\mc{L}\otimes \mc{N} | \cong \mc{O}_{\mc{E}}$ or $\mc{N} |_{\mc{E}} \cong \mc{O}_{\mc{E}}$.
  \item $\mc{C}$ is a smooth rational curve with $\Sigma^{\mc{C}} = \emptyset$ and $d_0 = d_{\infty} = 0$.
  \item $\mc{C}$ is irreducible, $g=1$, $\Sigma^{\mc{C}} = \emptyset$, $\mc{L}^k \cong \mc{O}_{\mc{C}}$, and $\mc{N} \cong \mc{L}^{-1}$.
  \end{enumerate}
\end{lem}

\begin{proof}
  Suppose that $\xi$ is unstable. By~\Cref{lem:stability}, we know that $\mc{C}$ has an unstable component $\mc{E}$. We first study the case where $g(\mc{E}) = 0$. First, suppose that $\mc{E}$ has one special point. If $\rho|_{\mc{E}} = 0$, then $\nu|_{\mc{E}}$ is nowhere zero, so $\mc{N}|_{\mc{E}} \cong \mc{O}_{\mc{E}}$. Then we know that $(\varphi|_{\mc{E}}, \mu|_{\mc{E}})$ defines a map $\mc{E} \to \P(\ba, 1^N)$ which must be constant, so $\mc{L}^k |_{\mc{E}} \cong \mc{O}_{\mc{E}}$.

  If $\rho |_{\mc{E}} \neq 0$, then $\deg \mc{L}|_{\mc{E}} < 0$, so $\varphi = 0$. Therefore $\mu|_{\mc{E}}$ is nowhere zero. Because it must be fixed by infinitely many automorphisms of $\mc{E}$, it is constant. This implies $\mc{L} \otimes \mc{N} |_{\mc{E}} \cong \mc{O}_{\mc{E}}$. Finally, we must show that $\rho$ is nowhere zero. Suppose that $\rho$ has a unique zero at $p_1$ and $\nu$ has a unique zero at $p_2 \neq p_1$ (this is the only possible case). Then either $p_1$ or $p_2$ is special. There is a $\C^{\times}_t$ acting on $\mc{E}$ by automorphisms $a_t$ and fixing $p_1, p_2$. Then if $p_1$ is special, we have isomorphisms
  \begin{align*}
    b_t &\colon \mc{L}|_{\mc{E}} \to a_t^* \mc{L}|_{\mc{E}} \\
    c_t &\colon \mc{N}|_{\mc{E}} \to a_t^* \mc{N}|_{\mc{E}}.
  \end{align*}
  Note that $b_t c_t = 1$. We then obtain the equations
  \begin{align*}
    b_t(p_2)^k t \rho(p_2) &= \rho(p_2) \\
    b_t(p_i)^{-1} \nu(p_i) &= \nu(p_i) \qquad i = 1,2.
  \end{align*}
  By nonvanishing of $\rho$ at $p_2$ and of $\nu$ at $p_1$, we see that $b_t(p_1) = 1$ and $b_t(p_2) = t^{1/k}$, so $\deg \mc{L}|_{\mc{E}} = -\frac{1}{k}$. But this implies that $\mc{L}^{-k}|_{\mc{E}} \otimes \omega_{\mc{E}}^{\log} \cong \mc{O}_{\mc{E}}$, so $\rho$ is nowhere vanishing. If $p_2$ is special, then the equations are
  \begin{align*}
    b_t(p_2)^k \rho(p_2) &= \rho(p_2) \\
    b_t(p_i)^{-1} \nu(p_i) &= \nu(p_i) \qquad i = 1,2.
  \end{align*}
  This gives us $b_t(p_1) = 1$ and $b_t(p_2)^{-k} = 1$, so $\deg \mc{L}|_{\mc{E}} = 0$, which is a contradiction to the earlier $\deg \mc{L}|_{\mc{E}} < 0$.

  Now suppose that $\mc{E}$ contains two special points. If $\rho|_{\mc{E}} = 0$, then the same argument as above shows that $\mc{L}^k|_{\mc{E}} \cong \mc{O}_{\mc{E}}$ and $\mc{N}|_{\mc{E}} \cong \mc{O}_{\mc{E}}$. If $\rho|_{\mc{E}} \neq 0$, then again suppose $\rho(p_1) = 0$ and $\nu(p_2) = 0$, where $p_1 \neq p_2$. Note that both $p_1, p_2$ must be special points (there must be infinitely many automorphisms fixing $p_1, p_2$, and the special points) and $p_1 \neq p_2$. Studying the action of $\C^{\times}_t$ and $b_t, c_t$ as before, we find the equations
  \begin{align*}
    b_t(p_2)^k \rho(p_2) &= \rho(p_2) \\
    b_t(p_i)^{-1} \nu(p_i) &= \nu(p_i) \qquad i = 1,2.
  \end{align*}
  Because $\rho(p_2) \neq 0$ and $\nu(p_1) \neq 0$, we obtain $b_t(p_1) = 1$ and $b_t(p_2)^k = 1$, so $\deg \mc{L}|_{\mc{E}} = 0$. But $\omega_{\mc{E}}^{\log} \cong \mc{O}_{\mc{E}}$, so $\mc{L}^k |_{\mc{E}} \cong \mc{O}_{\mc{E}}$. We again have $\varphi = 0$ because $\xi \in \mc{W}^{\mr{pre}-}$, so again $\mu$ is nowhere zero. This is possible only if $\mu$ is constant, so $\mc{L} \otimes \mc{N}|_{\mc{E}} \cong \mc{O}_{\mc{E}}$.

    If $g(\mc{E}) = 0$ and $\mc{E} \cap (\Sigma^{\mc{C}} \cup \mc{C}_{\mr{sing}}) = \emptyset$, then $\mc{C} = \mc{E} \cong \P^1$. If $\rho = 0$, then $\nu$ is nowhere zero, so $\mc{N} \cong \mc{O}_{\mc{C}}$. Then $(\varphi,\mu)$ defines a map $\mc{C} \to \P(\ba,1^N)$, which must be constant. This implies $\mc{L} \cong \mc{O}_{\mc{C}}$. If $\rho \neq 0$, then $\deg \mc{L} < 0$, so $\varphi = 0$. Then $\mu$ is nowhere zero, so it defines a map $\mc{C} \to \P^{N-1}$. This must be constant, so $\mc{L} \otimes \mc{N} \cong \mc{O}_{\mc{C}}$. If $\rho$ is not nowhere vanishing, as above suppose that $\rho(p_1) = 0$ and $\nu(p_2) = 0$, where $p_1 \neq p_2$. They must be fixed by infinitely many automorphisms of $\P^1$, so we may once again study the action of $\C^{\times}_t$ on $\mc{C} \cong \P^1$. We obtain the equations
  \begin{align*}
    b_t(p_2)^k t \rho(p_2) &= \rho(p_2) \\
    b_t(p_i)^{-1} \nu(p_i) &= \nu(p_i) \qquad i = 1,2.
  \end{align*}
  Because $\rho(p_2) \neq 0$ and $\nu(p_1) \neq 0$, we obtain $b_t(p_1) = 1$ and $b_t(p_2) = t^{1/k}$, so $\deg \mc{L} = -\frac{1}{k}$, which is impossible because $\mc{C} \cong \P^1$ is a scheme. Therefore, $\rho$ must be nowhere zero. This would imply that $k \leq 2$, which is a contradiction.

    Finally, we consider the case when $g(\mc{E}) = 1$. Then $\mc{E} \cap \Sigma^{\mc{C}} = \emptyset$ and $\mc{E} = \mc{C}$. If $\rho = 0$, then $\nu$ is nowhere zero, so $\mc{N} \cong \mc{O}_{\mc{C}}$. Again, $(\varphi,\mu)$ defines a map $\mc{C} \to \P(\ba,1^N)$ which must be constant, so $\mc{L} \otimes \mc{N} \cong \mc{O}_{\mc{C}}$. Otherwise, $\rho \neq 0$. Then $\varphi = 0$, so $\mu$ is nowhere vanishing and thus defines a map $\mc{C} \to \P^{N-1}$, which must be constant. This implies $\mc{L} \otimes \mc{N} \cong \mc{O}_{\mc{C}}$. If $\rho$ has a zero, then $\deg \mc{L} < 0$, so $\deg \mc{N} > 0$. Thus $\nu$ has a zero. Because $(\rho,\nu)$ must be nonvanishing, $(\rho = 0)$ and $(\nu = 0)$ are disjoint from each other. But any two points on a curve of arithmetic genus $1$ can only be fixed by finitely many automorphisms. Therefore, $\rho$ is nowhere vanishing, so $\mc{L}^k \cong \mc{O}_{\mc{C}}$. This puts us in case (4).

    We now prove the other direction. In case (1), we know $\deg \omega^{\log}_{\mc{C}}|_{\mc{E}} = -1$ because $\mc{E} \cap (\Sigma_{\mc{C}} \cap \mc{C}_{\mr{sing}})$ contains one point, say $p$. If $\rho|_{\mc{E}}$ is nowhere zero, then $\deg \mc{L}|_{\mc{E}} = -\frac{1}{k}$, so $p$ is stacky. Therefore $\mc{E} \cong \P(1,k)$, and because $\xi \in \mc{W}^-$, $\varphi|_{\mc{E}} = 0$. Following the same strategy as in the last paragraph of the proof of~\cite[Lemma 3.3]{mspfermat}, we see that $G = \G_a$ is a subgroup of the automorphism group of $\xi|_{\mc{E}}$. In the other cases, all of the fields are constant, so it is clear that $\xi$ is unstable.
\end{proof}

\subsection{Specialization of MSP fields}
\label{sub:specialization}

We will now study the specialization of MSP fields with the goal of proving that $\mc{W}^-$ is universally closed. From now, we will let $S$ be an affine smooth curve over an algebraically closed field with a closed point $s_0 \in S$ and let $S_* \coloneqq S \setminus s_0$. We will denote familes over $S_*$ by 
\[ \xi_* = (\mc{C}_*,\Sigma^{\mc{C}_*},\mc{L}_*,\mc{N}_*,\varphi_*,\rho_*,\mu_*,\nu_*) \in \mc{W}^{\pre}_{g,\gamma,\mbf{d}}. \]

Note as in~\cite[\href{https://stacks.math.columbia.edu/tag/0CL9}{Tag 0CL9}]{stacks-project} that in using the valuative criterion to prove properness, we need to take a finite base change $S' \to S$ ramified over $s_0$. If necessary, we can shrink $S$ and assume that there is an \'etale $S \to \A^1$ such that $s_0$ is the only point lying over $0 \in \A^1$, and then we can take $S' = S_a \coloneqq S \times_{\A^1,t \mapsto t^a} \A^1$. For notational convenience, when we do such a base change, we will replace $S'$ with $S$.

\begin{lem}
    Suppose $\rho = 0$. Then after a finite base change, $\xi_* \in \mc{W}^{-}(S_*)$ extends to a $\xi \in \mc{W}^-(S)$. 
\end{lem}

\begin{proof}
    Because $\rho_* = 0$, $\nu_*$ is nowhere vanishing, so $\mc{N}_* \cong \mc{O}_{\mc{C}_*}$. Therefore, $(\varphi_*, \mu_*)$ defines a stable map $f_* \colon \mc{C}_* \to \P(\ba,1^N)$ such that $(f_*)^* \mc{O}(1) \cong \mc{L}_*$. Because the moduli of stable maps is proper, after a finite base change we can extend $f_*$ to an $S$-family of stable maps $f \colon \mc{C} \to \P(\ba,1^N)$. Then we set $\mc{L} = f^* \mc{O}(1)$. By definition, $f$ is provided by a section $(\varphi,\mu) \in H^0\ab(\bigoplus \mc{L}^{a_j} \oplus \mc{L}^{\oplus N})$ extending $(\varphi_*,\mu_*)$. Finally, define $\mc{N} \cong \mc{O}_{\mc{C}}$ and $\nu$ to be the isomorphism $\mc{N} \cong \mc{O}_{\mc{C}}$ extending $\nu_*$ and define $\rho = 0$. This is the desired extension.
\end{proof}

\begin{lem}
    Suppose $\xi_* \in \mc{W}(S_*)$ satisfies $\varphi_* = \nu_* = 0$. Then after a finite base change, $\xi_*$ extends to $\xi \in \mc{W}(S)$.
\end{lem}

\begin{proof}
    Because $\varphi_* = \nu_* = 0$, $\rho_*$ and $\mu_*$ are nowhere vanishing. Therefore, we obtain a stable map
    \[ \mu_* \coloneqq \mc{C}_* \to \P^{N-1} \]
    and an $r$-spin structure, or in other words an $S_*$-family of stable spin maps as in~\cite{spingw}. By~\cite[Proposition 2.23]{spingw}, the moduli stack of spin stable maps is proper, so there is an extension of $(\mc{C}_*, \Sigma^{\mc{C}_*},\mc{L}_*,\mu_*)$ to a spin stable map $(\mc{C}, \Sigma^{\mc{C}},\mc{L},\mu \colon \mc{C} \to \P^{N-1})$ over $S$. Finally, we define $\mc{N} = \mu^* \mc{O}(1), \varphi = 0, \nu = 0$.
\end{proof}

We now work in the case where $\varphi_* = 0$ and $\rho_* \neq 0, \nu_* \neq 0$ and specialize to the case of smooth domain curves. Now $\mc{L}$ and $\mc{N}$ must have opposite monodromy at the marked points because $\mu_*$ is nowhere vanishing. But this implies that $\nu$ must vanish at the stacky points and that both $\mc{L}$ and $\mc{N}$ are representable, so $\rho$ cannot vanish at the stacky points. Therefore, we are in the situation of~\cite[\S 3.2]{nmsp}.

\begin{lem}\label{lem:extsmoothsource}
    Let $\xi_* \in \mc{W}^-(S_*)$ be as above. Then possibly after a finite base change, we can extend $\xi_*$ to $\xi \in \mc{W}^-(S)$.
\end{lem}

    

\begin{proof}
    The proof is exactly the same as in~\cite[\S 3.2]{nmsp}.
\end{proof}

We are now ready to prove the existence part of the valuative criterion for properness. After we prove that $\mc{W}^-_{g,\gamma,\mbf{d}}$ is finite type, this will imply that $\mc{W}^-_{g,\gamma,\mbf{d}}$ is universally closed. To avoid some technical difficulties, we will asssume that $S = \Spec R$ for $R$ a complete DVR with algebraically closed residue field using~\cite[Theorem 7.10]{champs}.\footnote{In~\cite{mspfermat}, $S$ is taken to be an arbitrary affine smooth curve, but there are some places where it is not clear to the author why certain gerbes and line bundles can be glued. To avoid this situation, we force everything to be trivial.}

\begin{prop}\label{prop:existence}
    Let $\xi_* \in \mc{W}^-(S_*)$. Then possibly after a finite base change, $\xi_*$ extends to $\xi \in \mc{W}^-(S)$. 
\end{prop}

\begin{proof}
    The proof is similar to that of~\cite[Proposition 3.23]{mspfermat} except in their setting, if $\nu_1$ is nowhere vanishing, then $\mc{L} \otimes \mc{N} \cong \mc{O}$. This is used in the case where $\rho \neq 0$ at a node.

    Possibly after a finite base change, we can assume that every connected component of the singular locus $(\mc{C}_*)_{\mr{sing}}$ is the image of a section of $\mc{C}_* \to S_*$. Let
    \[ \pi \colon \wt{\mc{C}}_* = \bigsqcup_\alpha \wt{\mc{C}}_{\alpha*} \to \mc{C}_* \]
    be the normalization of $\mc{C}_*$ where the $\wt{ \mc{C} }_{\alpha*}$ are the connected components. After a finite base change, we can asssume that $\wt{ \mc{C} }_{\alpha_*} \to S_*$ have connected fibers. Then at every field-valued point of $S_*$ every node of the corresponding fiber of $\mc{C}_*$ satisfies either $\varphi_* = 0$ or $\rho_* = 0$, so each
    \[ \xi_{\alpha*} \in \mc{W}^-_{g_{\alpha},\gamma_{\alpha},\mbf{d}_{\alpha}} \]
    for some $(g_{\alpha},\gamma_{\alpha},\mbf{d}_{\alpha})$. Applying~\Cref{lem:extsmoothsource}, after a finite base change $S_{\alpha} \to S$, we can extend $\xi_{\alpha*}$ to some $\xi_{\alpha} \in \mc{W}^-_{g_{\alpha},\gamma_{\alpha},\mbf{d}_{\alpha}}(S_{\alpha})$. We then choose a finite base change $S' \to S$ factoring through all $S_{\alpha}$, and replace $\xi_{\alpha}$ by $\xi_{\alpha} \times_{S_{\alpha}} S'$. Therefore, after a finite base change, we now have $\xi_{\alpha} \in \mc{W}^-_{g_{\alpha},\gamma_{\alpha},\mbf{d}_{\alpha}}(S)$ for all $\alpha$. Denote $\wt{\xi} = \bigsqcup_{\alpha} \xi_{\alpha}$.

    We now glue the $\xi_{\alpha}$ to obtain a stable $\xi$. Define $\wt{\mc{C}} = \bigsqcup_{\alpha} \mc{C}_{\alpha}$. Let $\Upsilon_* \subset \mc{C}_*$ be (the image of) a section of $(\mc{C}_*)_{\mr{sing}}$ and $\pi^{-1}(\Upsilon_*) \eqqcolon \Upsilon_{1*} \sqcup \Upsilon_{2*}$ be the corresponding markings in $\wt{\mc{C}}_*$. Possibly after a finite base change, $\Upsilon_{m*}$ extends to $\Upsilon_m$ for $m=1,2$. By assumption on the moduli spaces, there are trivializations of $\Upsilon_m$. Then we know that 
    \[ N_{\Upsilon_{1}/\wt{\mc{C}}} \cong N_{\Upsilon_2/\wt{\mc{C}}} \] 
    by construction (the normal bundle is simply a character of $\mu_a$, and the two agree on the generic point). Using~\cite[Definition 1.7]{orbdegeneration} when $\Upsilon_1,\Upsilon_2$ lie on different components and~\cite[Proposition A.1.1]{gwdmstack} when they lie on the same component, we obtain a gluing $\mc{C}$ of $\wt{\mc{C}}$ along $\Upsilon_1 \cong \Upsilon_2$, giving an $S$-family of (not necessarily connected) twisted curves. We can perform the gluing to all sections of $(\mc{C}_*)_{\mr{sing}}$, we obtain an $S$-family of twisted curves extending $\mc{C}_* \to S_*$. Denote this gluing morphism by
    \[ \beta \colon \wt{\mc{C}} \to \mc{C}. \]

    The next step is to glue the sheaves and fields. Without loss of generality, we may assume that $(\mc{C}_*)_{\mr{sing}}$ consists of a single $S_*$-section, so $\mc{C}$ is the gluing of $\wt{\mc{C}}$ along $\Upsilon_1 \cong \Upsilon_2$. Let 
    \[ \iota_m \colon \Upsilon \cong \Upsilon_m \to \wt{\mc{C}} \]
    be the tautological maps. The first case is when $\iota_1^* \wt{\varphi} \neq 0$. This implies both that $\iota_2^* \wt{\varphi} \neq 0$ and that $\wt{\rho} = 0$ in a neighborhood $\wt{\mc{U}} \subset \wt{\mc{C}}$ of $\Upsilon_1 \cup \Upsilon_2$ because we are working in $\mc{W}^-$. This implies $\wt{\nu}$ is nowhere zero on $\wt{\mc{U}}$, and induces an isomorphism $\wt{N}|_{\wt{\mc{U}}} \cong \mc{O}_{\wt{\mc{U}}}$. This induces $\iota_m^* \wt{\mc{N}} \cong \mc{O}_{\Upsilon}$ such that $\iota_m^* \wt{\nu} = 1$. Then
    \[ f_m \coloneqq (\iota_m^*\wt{\varphi}, \iota_m^* \wt{\mu}) \colon \Upsilon \to \P(\ba,1^N) \]
    satisfies $f_1 |_{\Upsilon_*} = f_2|_{\Upsilon_*}$, so $f_1 = f_2$. Therefore, there exists a unique isomorphism $\phi \colon \iota_1^* \wt{\mc{L}} \cong \iota_2^* \wt{\mc{L}}$ such that
    \[ \phi^* (\iota_2^* \wt{\varphi}, \iota_2^* \wt{\mu}) = (\iota_1^* \wt{\varphi}, \iota_1^* \wt{\mu}). \]
    Finally, $\wt{\mc{L}}^{-k} \otimes \omega_{\wt{\mc{C}}/S}^{\log}$ glues to $\mc{L}^{-k} \otimes \omega_{\mc{C}/S}^{\log}$. We have shown that
    \[ (\iota_1^* \wt{\varphi},\iota_1^* \wt{\rho}, \iota_1^* \wt{\mu}, \iota_1^* \wt{\nu}) = (\iota_2^* \wt{\varphi},\iota_2^* \wt{\rho}, \iota_2^* \wt{\mu}, \iota_2^* \wt{\nu}), \]
    so by~\cite[Corollary 3.22]{mspfermat}, we can glue $(\wt{\mc{L}},\wt{\mc{N}},\wt{\varphi},\wt{\rho},\wt{\mu},\wt{\nu})$ to $(\mc{L},\mc{N},\varphi,\rho,\mu,\nu)$.

    The second case is that $\iota_1^* \wt{\rho} \neq 0$. This forces $\iota_2^* \wt{\rho} \neq 0$. Because we are working in $\mc{W}^-$, $\wt{\varphi}|_{\wt{U}} = 0$ in a neighborhood $\wt{\mc{U}} \subset \wt{\mc{C}}$ of $\Upsilon_1 \cup \Upsilon_2$. Therefore,
    \[ f_m \coloneqq \iota_m^* \wt{\mu} \colon \Upsilon \to \P^{N-1} \]
    satisfies $f_1|_{\Upsilon_*} = f_2|_{\Upsilon_*}$, so there is a unique isomorphism 
    \[\phi \colon \iota_1^* ( \wt{\mc{L}} \otimes \wt{\mc{N}} ) \cong \iota_2^* ( \wt{\mc{L}} \otimes \wt{\mc{N}} ) \]
    such that $\phi^* \iota_1^* \wt{\mu} = \iota_2^* \wt{\mu}$. Next, we trivialize $\iota_m^* \omega_{\wt{\mc{C}}/S}^{\log}$ by $\frac{\d{x}}{x}$ (where $x$ is the fiberwise coordinate) and consider
    \[ \iota_m^* \wt{\rho} \in H^0(\iota_m^* \wt{\mc{L}}^{-k}). \]
    Using again that $S$ is the spectrum of a local ring, we trivialize $\iota_m^* (\wt{\mc{L}} \otimes \wt{\mc{N}}) \cong \mc{O}_{\Upsilon}$ in a way compatible with $\phi$. Therefore, $\iota_m^* \wt{\mc{N}} \cong \iota_m^* \wt{\mc{L}}^{-1}$, so
    \[ h_m \coloneqq (\iota_m^* \wt{\rho}, \iota_m^* \wt{\nu}) \colon \Upsilon_m \to \P(1,k) \]
    satisfies $h_1 |_{\Upsilon_*} = h_2|_{\Upsilon_*}$. Therefore $h_1 = h_2$ and there exists a unique isomorphism $\phi' \colon \iota_1^* \wt{\mc{L}}^{-1} \cong \iota_2^* \wt{\mc{L}}^{-1}$ such that
    \[ \phi^* (\iota_2^* \wt{\rho}, \iota_2^* \wt{\nu}) = (\iota_1^* \wt{\rho}, \iota_2^* \wt{\nu}). \]
    We have again shown that
    \[ (\iota_1^* \wt{\varphi},\iota_1^* \wt{\rho}, \iota_1^* \wt{\mu}, \iota_1^* \wt{\nu}) = (\iota_2^* \wt{\varphi},\iota_2^* \wt{\rho}, \iota_2^* \wt{\mu}, \iota_2^* \wt{\nu}), \]
    so applying~\cite[Corollary 3.22]{mspfermat} as above, we can glue $(\wt{\mc{L}},\wt{\mc{N}},\wt{\varphi},\wt{\rho},\wt{\mu},\wt{\nu})$ to $(\mc{L},\mc{N},\varphi,\rho,\mu,\nu)$. The last step is to set $\Sigma^{\mc{C}} = \beta(\Sigma^{\wt{\mc{C}}} \setminus ( \Upsilon_1 \cup \Upsilon_2 ))$.
\end{proof}

\subsection{Valuative criterion for separatedness}

Again, let $S$ be an affine smooth curve over an algebraically closed field and $s_0 \in S$ be a closed point.

\begin{lem}\label{lem:smuniqueness}
    Let $\xi, \xi' \in \mc{W}(S)$ be such that $\xi_* \cong \xi'_*$ and $\mc{C}_*$ is smooth. Then $\xi \cong \xi'$.
\end{lem}

\begin{proof}
    Let $C$ (resp. $C'$) be the coarse moduli space of $\mc{C}$ (resp. $\mc{C}'$). Let $\pi \colon X \to C$ (resp. $\pi' \colon X' \to C'$) be the minimal desingularization. Because $C,C'$ have only type $A$ singularities, $\pi, \pi'$ are contractions of chains of $(-2)$-curves.

    Let $f \colon X \dashrightarrow X'$ be the birational map induced by $\xi_* \cong \xi'_*$ and let $U_0 \subset X$ be the largest open subset over which $f$ is defined. If $U_0 \neq X$, then $X \setminus U_0$ is discrete. Let $X_1$ be the blowup of $X$ at $X \setminus U_0$. Inductively, suppose $X_m \to X$ is a successive blowup and let $U_m \subset X_m$ be the largest open subset over which $f \colon X_m \dashrightarrow X'$ is defined. Then define $X_{m+1}$ to be the blowup of $X_m$ at $X_m \setminus U_m$. After finitely many steps, there exists $\bar{m}$ such that $Y \coloneqq X_{\bar{m}} = U_{\bar{m}}$. Denote the induced projection and birational morphism by
    \[ \bar{\pi} \colon Y \to X, \qquad \bar{f} \colon Y \to X'. \]

    Now write $E \subset Y$ for the exceptional divisor of $\bar{\pi}$. By construction, we can write
    \[ E = \sum_{m \geq 1} E_m, \]
    where $E_m \subset E$ is the strict transform of the exceptional divisor of $X_m \to X_{m-1}$, where $X_0 \coloneqq X$. Let $E' \subset Y$ be the exceptional divisor of $\bar{f}$ and $V = Y \setminus E'$. By construction, $E$ and $E'$ share no common irreducible components. Also, let
    \[ Y_0 = \bigcup_j D_j \]
    be the decomposition of the central fiber $Y_0$ of $Y$. Using the formula for the dualizing sheaf of a blowup, we have
    \begin{equation}\label{eqn:blowup}
        \omega_{Y/S}^{\log} \cong \bar{\pi}^* \omega_{X/S}^{\log}\ab(\sum_i i E_i). 
    \end{equation}
    Note that $\mc{L}^k, \mc{N}^k$ (resp. $\mc{L}'^k, \mc{N}'^k$) descend to $C$ (resp $C'$). Then define
    \[ L \coloneqq \pi^* \mc{L}^k, \qquad M \coloneqq \pi^* \mc{N}^k \]
    and similarly for $L'$ and $M'$. Also let $h_j$, $u_{l}$, and $v$ be the pullbacks of $\varphi_j^{b_j}$, $\mu_{l}^k$, and $\nu^k$, respectively, which are sections of $L, L \otimes M$, and $M$, respectively. Denote $h_j', u_{l}'$, and $v'$ simlarly. We will also view $\rho \in H^0(X, L^{-1} \otimes \omega_{X/S}^{\log})$ and $\rho'$ similarly. Because $\xi_* = \xi'_*$, there exist integers $a_i, b_i$ such that
    \begin{equation}\label{eqn:pullbacks}
        \bar{f}^* L' \cong \bar{\pi}^* L\ab(\sum a_i D_i), \qquad \bar{f}^* L' \cong \bar{\pi}^* L\ab(\sum b_i D_i). 
    \end{equation}

    First, suppose that $E \neq \emptyset$. Then let $D_j \subset E$ be an irreducible component and $x = \bar{\pi}(D_j) \in X$. If $x$ is a singular point of $X_0$, then $D_j$ must have multiplicity $\geq 2$ as a component of $Y_0$ and thus $Y_0$ is not reduced. But $D_j \cap V$ is $1$-dimensional and $\bar{f}|_V$ is an isomorphism onto its image, so this would imply that $X_0'$ is not reduced. But $X_0'$ is reduced by construction, so $x$ must be a smooth point of $X_0$. We know that $\bar{\pi}^{-1}(x)$ is a tree of rational curves. We may assume that $D_1 + \cdots + D_m$ form a maximal chain of rational curves in $\bar{\pi}^{-1}(x)$ possibly after reindexing the components of $Y_0$. We may also assume that $D_i \subseteq E_i$ for all $i$. By construction, $D_m$ and $\bar{f}(D_m)$ are $(-1)$-curves, so $\mc{D}' \coloneqq \pi' (\bar{f}(D_j))$ is a rational curve. Finally, let $z \in D_m$ be a general point and $y = \bar{f}(z)$.

    \begin{claim}\label{claim:viszero}
        In the situation as above, $v(x) = 0$.
    \end{claim}

    \begin{proof}
        Suppose $u_2(x) \neq 0$. Then suppose that $u_{l}(x) \neq 0$ for some $l$. Then $\bar{\pi}^* v |_{D_m} \neq 0$ and $\bar{\pi}^* u_{l}|_{D_m} \neq 0$, while $\bar{f}^* v'$ and $\bar{f}^* u_{l}'$ could vanish along $D_m$. Because $\bar{\pi}^* u_{l} |_{V_*} = \bar{f}^* u_{l}'|_{V_*}$ and $\bar{\pi}^* v |_{V_*} = \bar{f}^* v'|_{V_*}$, we see that $b_m, a_m + b_m \geq 0$.

        If $a_m + b_m > 0$, then $\bar{f}^*u_{l}'|_{D_m} = 0$, so $u_{l}'(y) = 0$. But then $v'(y) \neq 0$, so $b_m \leq 0$. This argument works for all $l$ such that $u_{l}(x) \neq 0$, so we see that $u'(y) = 0$. This implies that there exists $j$ such that $h_j'(y) \neq 0$, which forces $a_m \leq 0$, a contradiction. Therefore $a_m + b_m = 0$.

        If $b_m > 0$, then $v'(y) = 0$ and $a_m < 0$. Using~\Cref{eqn:blowup} and~\Cref{eqn:pullbacks}, there exists a dense open subset $U \subset D_i$ such that
        \begin{equation}\label{eqn:logcansingledivisor}
            \bar{f}^*(L'^{-1} \otimes \omega_{X'/S}^{\log})|_U \cong \bar{\pi}^*(L^{-1} \otimes \omega^{\log}_{X/S})((i-a_i)D_i) |_{U}
        \end{equation}
        for any $i \leq m$. Setting $i=m$, we obtain $\rho'(y) = 0$. But $v'(y) = 0$, which is a contradiction. Therefore, $a_m = b_m = 0$.

        In the case where $u(x) = 0$, there exists $j$ such that $h_j(x) \neq 0$. Also, $v(x) \neq 0$. As above, we conclude that $a_m, b_m \geq 0$. If $b_m > 0$, then $\bar{f}^* v'|_{D_m} = 0$, so $\bar{f}^* \rho'|_{D_m} \neq 0$. Applying~\Cref{eqn:logcansingledivisor}, we see that $m - a_m \leq 0$, so $a_m \geq m > 0$. This implies that $a_m + b_m > 0$, so $\bar{f}^* u'|_{D_m} = 0$, which is a contradiction because $(u',v')$ must be nowhere zero. Therefore, $b_m = 0$. If $a_m > 0$, then $\bar{f}^* h'|_{D_m} = 0$. But because $b_m = 0$ and $\bar{\pi}^*u|_{D_m} = 0$, we also have $\bar{f}^* u|_{D_m} = 0$ (otherwise $b_m < 0$). This contradicts the nonvanishing of $(h',u')$, so $a_m = 0$ as well.
        
        Finally, because $a_m = 0$, we know that $\bar{f}^*\rho'|_{D_m} = 0$. Descending to $\mc{D}'$, we see that $\rho'|_{\mc{D}'} = 0$ and thus $\mc{N}'|_{\mc{D}'} \cong \mc{O}_{\mc{D}'}$. Then
        \[ (\varphi', \mu') \colon \mc{D}' \to \P(\ba,1^N) \]
        is a morphism. But on a dense open subset $U$ of $D_m$ away from the nodes of $Y_0$, we have
        \[ \bar{f}^* L'|_U \cong \bar{\pi}^*L |_U, \qquad \bar{f}^* M'|_U \cong \bar{\pi}^*M |_U \]
        because $a_m = b_m = 0$ and thus 
        \[ (\bar{f}^* h', \bar{f}^*u') = (\bar{\pi}^* h, \bar{\pi}^* u). \]
        But now $\bar{\pi}(D_m)$ is a point, so $(\bar{\pi}^* h, \bar{\pi}^* u)$ is a constant. This implies that $(\varphi',\mu')|_{\mc{D}'}$ is constant, so $\deg \mc{L} = \deg \mc{N} = 0$. Because $\mc{D}'$ contains at most two special points (one node and one marked point), by~\Cref{lem:stabilitydegeneracy}, $\xi_0'$ is unstable.
    \end{proof}

    \begin{claim}
        We have $a_m = m$, $a_i \geq a_{i+1}-2$ for $i < m$, and $a_i + b_i = 0$ for all $i \leq m$.
    \end{claim}
    
    \begin{proof}
        First, by the previous claim, we know that $v(x) = 0$. This implies that $\rho(x) \neq 0$, so $\bar{\pi}^* \rho|_{D_i} \neq 0$. Using~\Cref{eqn:logcansingledivisor}, we see that $i - a_i \geq 0$ for all $i$. Then $u(x) \neq 0$, so $a_i + b_i \geq 0$.

        Now if $a_i + b_i > 0$, this forces $\bar{f}^* u'|_{D_i} = 0$. Then $\bar{f}^* h'|_{D_i} \neq 0$, so $a_i \leq 0$. Because $a_i + b_i > 0$, we must have $b_i > 0$, so $\bar{f}^* u_2'|_{D_i} = 0$, which implies that $(u'(y), v'(y)) = 0$, a contradiction. Therefore, $a_i + b_i = 0$.

        If $a_m \neq m$, then we know that $a_m < m$. Therefore, $\bar{f}^* \rho'|_{D_m} = 0$, so $\bar{f}^*v'|_{D_m}$ is nowhere vanishing. By the previous claim, $v(x) = 0$, so $\bar{\pi}^* v|_{D_m} = 0$ and thus $b_m < 0$. Because $a_m + b_m = 0$, we have $a_m > 0$. This implies that $\bar{f}^* h'|_{D_m} = 0$, so $\varphi'|_{\mc{D}'} = 0$. Adding that $\rho'|_{\mc{D}'} = 0$, we know that $\nu'|_{\mc{D}'}$ is nowhere vanishing. Because $a_m + b_m = 0$, we have $\bar{f}^* u|_{U} = \bar{\pi}^* u|_{U}$ for an open subset $U \subset D_m$ avoiding the nodes of $Y_0$. But $D_m$ is contracted by $\bar{\pi}$, so $\bar{f}^* u'|_{D_m}$ is constant. This implies that $\mu'|_{\mc{D'}}$ is constant. Because $\mc{D}'$ contains at most two special points, $\xi_0'$ is unstable by~\Cref{lem:stabilitydegeneracy}. Therefore, $a_m = m$.

        The final step is to prove that $a_i \leq a_{i+1}-2$. Let 
        \[ \Lambda \coloneqq \ab\{ 1 \leq i < m \mid a_i > a_{i+1}-2\}. \]
        If $\Lambda \neq \emptyset$, let $i$ be the largest element in $\Lambda$. If $i = m-1$, then $a_{m-1} \geq a_m - 1 =  m-1$. But we already know that $a_{m-1} \leq m-1$, so $a_{m-1} = m-1$. We conclude that
        \[ \bar{f}^*(L' \otimes \omega^{\log}_{X'/S}) \cong \bar{\pi}^* (L \otimes \omega^{\log}_{X/S}), \]
        so $\bar{f}^* \rho'|_{D_m} = \bar{\pi}^* \rho|_{D_m}$. Because $v(x) = 0$, $\rho'(x) \neq 0$, so $\bar{f}^* \rho'|_{D_m}$ is nowhere vanishing. Because $a_m + b_m = 0$, $\mu'|_{\mc{D}'}$ is constant as above. This implies $(\mc{L}' \otimes \mc{N}')|_{\mc{D}'} \cong \mc{O}_{\mc{D}'}$. As above, $\mc{D}'$ can contain either one or two special points of $\mc{C}_0'$, and in either case, $\xi'_0$ is unstable by~\Cref{lem:stabilitydegeneracy}. Thus $i < k-1$.

        Now we have $(i+1) - a_{i+1} > 0$, which implies $\bar{f}^* \rho'|_{D_{i+1}} = 0$ and thus $\bar{f}^* v'|_{D_{i+1}}$ is nowhere vanishing, so $\bar{f}^*M'|_{D_{i+1}} \cong \mc{O}_{D_{i+1}}$. Because $D_{i+1} \subseteq E$, $\bar{\pi}^*(L \otimes M) \cong \mc{O}_{D_{i+1}}$. Because $a_{i+1} + b_{i+1} = 0$, we obtain
        \[ \bar{f}^* (L' \otimes M') \cong \bar{\pi}^* (L \otimes M) \cong \mc{O}_{D_{i+1}}. \]
        Therefore, $\bar{f}^*L' |_{D_{i+1}} \cong \mc{O}_{D_{i+1}}$. If $D_i, D_{i+2}, D_{k_2}, \ldots, D_{k_m} \subseteq E$ be the irreducible components that intersect $D_{i+1}$. We know $a_{i+1} \leq a_{i+2}-2$. We may assume without loss of generality that $a_{i+1} \leq a_{k_s} - 2$ for all $s = 2,\ldots,m$. Using~\Cref{eqn:pullbacks} and $\bar{f}^* L' \cong \bar{\pi}^* L \cong \mc{O}_{D_{i+1}}$, we obtain
        \begin{align*}
            \mc{O}_{D_{i+1}} &\cong \mc{O}_{D_{i+1}}\ab(a_i D_i|_{D_{i+1}} + a_{i+1} D_{i+1}|_{D_{i+1}} + a_{i+2}D_{i+2}|_{D_{i+1}} + \sum_s a_{k_s} D_{k_s}|_{D_{i+1}}) \\
            &\cong \mc{O}_{D_{i+1}}\ab(a_i + a_{i+2} + \sum_s a_{k_s} - (2+s)a_{i+1}).
        \end{align*}
        This implies that
        \[ (a_i - a_{i+1}) + (a_{i+2} - a_{i+1}) + \sum_s (a_{k_s} - a_{i+1}) = 0, \]
        and combined with the fact that every term except the first is at least $2$, we obtain $a_i - a_{i+1} \leq -2$, a contradiction. Therefore, $\Lambda = \emptyset$, so $a_i \leq a_{i+1} - 2$ for all $i$.
    \end{proof}
    
    We continue with our chain of maximal curves $D_1, \ldots, D_m \subset E$. If $m \geq 2$, then $a_1 \leq a_2 - 2 \leq 0$. This implies $\bar{f}^* \rho'|_{D_1} = 0$, so $\bar{f}^* v'|_{D_1}$ is nowhere vanishing. By~\Cref{claim:viszero}, we know $\bar{\pi}^* v|_{D_1} = 0$, so $b_1 < 0$. But we know $a_1 = -b_1 > 0$, which is a contradiction. This implies $m=1$.

    Now let $D_j \subset Y_0$ not be contained in $E$. Suppose there exists $D_i \subset E$ such that $D_i \cap D_j \neq \emptyset$. If $\bar{f}$ does not contract $D_j$, then $a_j = b_j = 0$ and $D_j \cap V$ is dense in $D_j$. Let $x = \bar{\pi}(D_i)$. Then $v(x) = 0$, so $\bar{\pi}^* \rho|_{D_i}$ is nowhere vanishing. Because $a_i = 1$, $\bar{f}^* \rho'|_{D_i}$ is also nowhere vanishing. Combined with $a_j = 0$, this implies that $\deg(\bar{f}^* L'|_{D_i}) = -1$. Therefore, $\bar{f}^* h'|_{D_i} = \bar{\pi}^*h|_{D_i} = 0$. Because $a_i + b_i = 0$ and $D_i$ is contracted by $\bar{\pi}$, we also know $\bar{f}^* u|_{D_i} = \bar{\pi}^* u|_{D_i}$ is constant. But then $\xi'_0$ is unstable by~\Cref{lem:stabilitydegeneracy}.

    Because $Y$ resolves both $f$ and $f^{-1}$, everything we have done so far also applies to $E'$. Therfore, $E$ and $E'$ are both disjoint unions of $(-1)$-curves and $E \cup E' = Y_0$. This is possible only if $E \cong E' \cong \P^1$ and $Y$ is the blowup of $X = S \times \P^1$ at a single point in $X_0$. This would imply that either $\xi_0$ or $\xi'_0$ is unstable, which is a contradiction. Therefore, $E$ and $E'$ are both empty, so $f$ and $f^{-1}$ are birational morphisms. We have proven that $f$ is an isomorphism. Because any $D_i \subset Y_0$ satisfies $a_i = b_i = 0$, we see that
    \begin{equation}\label{eqn:cmsiso} 
        f^* L' \cong L, \qquad f^* M' \cong M, \qquad f^*(h',\rho',u',v') = (h,\rho,u,v). 
    \end{equation}

    Our goal now is to descend to an isomorphism $C \cong C'$. We note that if $D_i \subset X$ is contracted by $\pi \colon X \to C$, then $D_i$ is a $(-2)$-curve and $L|_{D_i} \cong M|_{D_i} \cong \mc{O}_{D_i}$. But if $D_i$ is a $(-2)$ curve with $L|_{D_i} \cong M|_{D_i} \cong \mc{O}_{D_i}$, $\xi|_{D_i}$ is unstable, so it must be contracted. Therefore, we obtain an isomorphism $\bar{\phi} \colon C \cong C'$.

    Let $\Delta$ (resp $\Delta'$) be the set of singular points of $\mc{C}_0$ (resp. $\mc{C}_0'$) and $\mr{pr} \colon \mc{C} \to C$ (resp. $\mr{pr}' \colon \mc{C'} \to C'$) be the coarse moduli space morphism. The isomorphism $\bar{\phi}$ and those of~\Cref{eqn:cmsiso} induce an isomorphism
    \[ \bar{\phi}^* \on{pr}'_* (\mc{L}'^k, \mc{N}'^k, \varphi_j'^{b_j}, \rho', \mu_{l}'^k, \nu'^k) \cong \on{pr}_* (\mc{L}^k, \mc{N}^k, \varphi_j^{b_j}, \rho,\mu_{l}^k, \nu^k). \]
    They also induce an isomorphism $\phi \colon \mc{C} \setminus \Delta \cong \mc{C'} \setminus \Delta'$ and isomorphisms
    \[ \phi^* (\mc{L'}, \mc{N'}, \varphi', \rho', \mu', \nu') \cong (\mc{L}, \mc{N}, \varphi, \rho, \mu,\nu) \]
    extending $\xi_* \cong \xi'_*$.

    The last step is to extend the isomorphism over $\Delta$. First, let $p \in \Delta$ and $p' \in \Delta'$ be the corresponding points. Choose an open subset $p \in \mc{U} \subset \mc{C}$ such that $\mc{U} \cap \Delta = p$. Then let $\mc{U'} \subset \mc{C'}$ be such that $\phi(\mc{U} \setminus p) = \mc{U'} \setminus p'$.

    Suppose that $\nu(p) \neq 0$. Then $\nu'(p') \neq 0$ as well. Therefore, after shrinking $\mc{U}$ if necessary, $\nu$ is nowhere vanishing on $\mc{U}$, providing an isomorphism $\mc{N} \cong \mc{O}_{\mc{U}}$. Similarly, after shrinking $\mc{U'}$, $\nu'$ is nowhere vanishing on $\mc{U'}$, so $\mc{N'} \cong \mc{O}_{\mc{U}}$. If $\nu$ and $\nu'$ are nowhere vanishing on all of $\mc{C}$ and $\mc{C}'$, respectively, then $(\varphi, \mu)$ and $(\varphi', \mu')$ define stable maps to $\P(\ba,1^N)$ (here, we are using the fact that $\mc{C}_*$ is smooth), and the extension is unique. Otherwise, if $\mu_{l}(p) \neq 0$ for all $l$, $\mu'_{l}(p') \neq 0$ for all $l$. This implies that $p$ and $p'$ are scheme points, so the uniqueness of extensions is straightforward.

    The other case is when $\mu(p) = 0$. If $\mu = 0$ everywhere, then $\nu$ is nowhere vanishing, which was discussed previously. Therefore, we may assume that there exists $l$ such that $(\mu_{l} = 0)$ and $(\mu'_{l} = 0)$ are divisors in $\mc{C}$ and $\mc{C'}$, respectively. By the same argument as in the end of the proof of~\cite[Lemma 4.1]{mspfermat}, the extension is unique.

    The other case is where $\nu(p) = \nu'(p') = 0$. When $\nu = 0$ everywhere, then $\nu' = 0$ everywhere as well. Then $\rho$ and $\rho'$ are nowhere vanishing, and similarly with $\mu$ and $\mu'$, which reduces to the case of spin stable maps to $\P^{N-1}$. By the results of~\cite{spingw}, the extension is unique.

    The case when $(\nu = 0)$ and $(\nu' = 0)$ are divisors in $\mc{C}$ and $\mc{C'}$, respectively, are treated as in~\cite[Lemma 4.1]{mspfermat}.
\end{proof}

\begin{prop}
    The result of~\Cref{lem:smuniqueness} holds without assuming that $\mc{C}_*$ is smooth.
\end{prop}

\begin{proof}
    The proof is exactly the same as in~\cite[Proposition 4.4]{mspfermat}.
\end{proof}

Using~\cite[Theorem 7.8]{champs}, we conclude that

\begin{prop}\label{prop:separated}
    The stack $\mc{W}_{g,\gamma,\mbf{d}}$ is separated.
\end{prop}

\subsection{$\mc{W}^-_{g,\gamma,\mbf{d}}$ is of finite type}

We now consider the fixed locus $(\mc{W}^-_{g,\gamma,\mbf{d}})^{\T}$ and study its $\C$-points. Let
\[ \xi = (\mc{C}, \Sigma^{\mc{C}}, \mc{L}, \mc{N}, \varphi,\rho,\mu,\nu) \in (\mc{W}^-_{g,\gamma,\mbf{d}})^{\T}(\C) \]
and $\mc{E}$ be an irreducible component of $\mc{C}$. Then either $\mc{E}$ carries a nontrivial action of $\T$ (in which case $g(\mc{E}) = 0$), or $\mu|_{\mc{E}} = 0$, or $\varphi|_{\mc{E}} = \rho|_{\mc{E}} = 0$, or $\nu|_{\mc{E}} = 0$. Let $\mc{C}_0$ be the one-dimensional part of $\mc{C} \cap (\mu = 0)$, $\mc{C}_1$ be the one-dimensional part of $\mc{C} \cap (\varphi = \rho = 0)$, and $\mc{C}_{\infty}$ be the one-dimensional part of $\mc{C} \cap (\nu = 0)$, all with the reduced structure. Then let $\mc{C}_{01}$ be the union of irreducible components in $\mc{C} \cap (\rho = 0)$ that are not in $\mc{C}_0 \cup \mc{C}_1 \cup \mc{C}_{\infty}$ and $\mc{C}_{1\infty}$ be the union of irreducible components in $\mc{C} \cap (\varphi = 0)$ that are not in $\mc{C}_0 \cup \mc{C}_{1} \cup \mc{C}_{\infty}$. By construction, we have:

\begin{lem}\leavevmode
    \begin{enumerate}
        \item $\mc{C}_0, \mc{C}_{01}, \mc{C}_1, \mc{C}_{1\infty}, \mc{C}_{\infty}$ contain no common irreducible components;
        \item The curves $\mc{C}_0, \mc{C}_1, \mc{C}_{\infty}$ are pairwise disjoint;
        \item $\mc{C} = \mc{C}_0 \cup \mc{C}_{01} \cup \mc{C}_1 \cup \mc{C}_{1\infty} \cup \mc{C}_{\infty}$.
    \end{enumerate}
\end{lem}


Now let $\Upsilon_{\xi}$ be the dual graph of $(\mc{C}, \Sigma^{\mc{C}})$ as in~\cite[\S 4.2]{mspfermat}. For $a = 0$, $01$ ,$1$, $1\infty$, or $\infty$, denote 
\[ V(\Upsilon_{\xi})_a = \ab\{v \in V(\Upsilon_{\xi}) \mid \mc{C}_v \subseteq \mc{C}_a\}. \]

\begin{lem}
    The set $\Pi \coloneqq \ab\{\Upsilon_{\xi} | \xi \in (\mc{W}^-_{g,\gamma,\mbf{d}})^{\T}(\C)\}$ is finite.
\end{lem}

\begin{proof}
    Note that the curve $(\mc{C}, \Sigma^{\mc{C}})$ may not be stable, so we use the information of $\mc{L}$ and $\mc{N}$ to form a (semi)stable $\wt{\Upsilon}_{\xi}$ by adding legs to $\Upsilon_{\xi}$.

    First, for every $v \in V(\Upsilon_{\xi})$, we add $3k \deg (\mc{L} \otimes \mc{N})|_{\mc{C}_v}$ auxillary legs to $v$. Here, note that $\deg (\mc{L} \otimes \mc{N})|_{\mc{C}_v} \geq 0$ because:
    \begin{enumerate}
        \item If $\mc{C}_v \subseteq \mc{C}_0$, then $\varphi$ and $\nu$ are nonzero;
        \item If $\mc{C}_v \subseteq \mc{C}_{01} \cup \mc{C}_1 \cup \mc{C}_{1\infty} \cup \mc{C}_{\infty}$, then $\mu$ is nonzero.
    \end{enumerate}
    The total number of new legs added is at most $3kd_0$, and a routine check shows that every vertex $v \notin V(\Upsilon_{\xi})_{1\infty}$ becomes stable.

    Now we define
    \[ E_{\infty} = \ab\{e \in E(\Upsilon_{\xi}) \mid e \in E_v \text{ for some } v \in V(\Gamma_{\xi})_{\infty}\}. \]
    For any $v \in V(\Upsilon_{\xi})_{1\infty}$, define
    \[ \delta(v) = r \deg \mc{N}|_{\mc{C}_v} - \ab\vert E_v \cap E_{\infty}\vert . \]
    Here, we note that $\deg \mc{N}|_{\mc{C}_v} > 0$ because otherwise $\xi|_{\mc{C}_v}$ is unstable, so $\delta(v) \in \Z_{\geq 0}$ for all $v \in V(\Upsilon_{\xi})$. To each such $v$, we add $2 \delta(v)$ auxillary legs to $v$.

    Let $\mc{C}_{\infty}^1, \ldots, \mc{C}_{\infty}^s$ be the connected components of $\mc{C}_{\infty}$ and $\ell_m$ be the number of markings on $\mc{C}_{\infty}^m$. Because $\nu|_{\mc{C}_{\infty}} = 0$, we know $\rho|_{\mc{C}_{\infty}}$ is nowhere vanishing. Because $\deg (\mc{L}\otimes \mc{N})|_{\mc{C}_{\infty}^m} \geq 0$, we see that
    \begin{align*}
        0 &= -k \deg \mc{L}|_{\mc{C}_{\infty}^m} + \deg \omega_{\mc{C}_{\infty}^m}^{\log} \\
        &\leq k \deg \mc{N}|_{\mc{C}_{\infty}^m} + 2g(\mc{C}_{\infty}^m) - 2 + \ab\vert \mc{C}_{1\infty} \cap \mc{C}_{\infty}^m\vert + \ell_m.
    \end{align*}
    Next, because $\sum_{v \in V(\Upsilon_{\xi})_{1\infty}} \ab|E_v \cap E_{\infty}| = \sum_{m=1}^s \ab|\mc{C}_{\infty}^m \cap \mc{C}_{1\infty}|$ and $\mc{N}|_{\mc{C}_v} \cong \mc{O}_{\mc{C}_v}$ unless $v \in V(\Upsilon_{\xi})_{1\infty} \cup V(\Upsilon_{\xi})_{\infty}$, we obtain
    \begin{align*}
        d_{\infty} &= \deg \mc{N} \\
        &\geq \frac{1}{k} \sum_{m=1}^s \ab( 2-2g(\mc{C}_{\infty}^m) - \ab|\mc{C}_{\infty}^m \cap \mc{C}_{1\infty}| - \ell_m) + \sum_{v \in V(\Upsilon_{\xi})_{1\infty}} \deg \mc{N}|_{\mc{C}_v} \\
        &\geq \frac{2s}{k} - \frac{2}{k}\sum_{m=1}^s g(\mc{C}_{\infty}^m) - \frac{1}{k} \ell + \frac{1}{k} \sum_{v \in V(\Upsilon_{\xi})_{1\infty}} \delta(v).
    \end{align*}
    Rearranging, we obtain
    \[ 2kd_{\infty} + 4g + 2\ell \geq 4s + 2\sum_{v \in V(\Upsilon_{\xi})_{1\infty}} \delta(v), \]
    and therefore the total number of auxillary legs and the number of connected components of $\mc{C}_{\infty}$ are bounded by $2kd_{\infty} + 4g + 2\ell$.

    Let $\wt{\Upsilon}_{\xi}$ be the resulting graph. If $v \in V(\wt{\Upsilon}_{\xi})$ is unstable, then $v \in V(\Upsilon_{\xi})_{1\infty}$. In addition, we must have $\delta(v) = 0$ and $\ab|E_v| = 1$ or $\ab|E_v| = 2$. If $\ab|E_v| = 1$, then because $\deg \mc{N}|_{\mc{C}_v} > 0$, we must have $\deg \mc{N}|_{\mc{C}_v} = \frac{1}{k}$. As in~\cite[Lemma 4.4]{nmsp}, we know that $(\mc{L} \otimes \mc{N})|_{\mc{C}_v} \cong \mc{O}_{\mc{C}_v}$, and therefore $\deg \mc{L}|_{\mc{C}_v} = -\frac{1}{k}$. This implies that $\rho$ is nowhere vanishing, which implies that $\xi$ is unstable. Therefore, $\ab|E_v| = 2$ and $\deg \mc{N}|_{\mc{C}_v} = \frac{1}{k}$, so $v$ is a strictly semistable vertex of $\wt{\Upsilon}_{\xi}$. As in the proof of~\cite[Proposition 4.6]{mspfermat}, there cannot be a chain of semistable vertices of length greater than $2$.

    Let $\wt{\Upsilon}_{\xi}^{\mr{st}}$ be the stabilization of $\wt{\Upsilon}_{\xi}^{\mr{st}}$. Because the genus is $g$ and the number of marked points is bounded, the set
    \[ \Pi^{\mr{st}} = \ab\{ \wt{\Upsilon}_{\xi}^{\mr{st}}\} \]
    is finite. Because we contracted chains of at most $2$ semistable vertices, $\Pi$ is finite.
\end{proof}

Because the (spin) Gromov-Witten moduli spaces are finite type, we obtain
\begin{lem}
    The fixed locus $(\mc{W}^-_{g,\gamma,\mbf{d}})$ is of finite type.
\end{lem}

\begin{prop}\label{prop:finitetype}
    $\mc{W}^-_{g,\gamma,\mbf{d}}$ is of finite type.
\end{prop}

\begin{proof}
    The proof is the same as that of~\cite[Proposition 3.9]{nmsp}.
\end{proof}

\subsection{Properness of the degeneracy locus}

\begin{thm}
    The stack $\mc{W}^-_{g,\gamma,\mbf{d}}$ is proper.
\end{thm}

\begin{proof}
    Because $\mc{W}^-_{g,\gamma,\mbf{d}}$ is separated and finite type by~\Cref{prop:separated} and~\Cref{prop:finitetype}, it suffices to check the existence part of the valuative criterion. This is exactly the statement of~\Cref{prop:existence}.
\end{proof}

\printbibliography

\end{document}